\documentclass[11pt, draftclsnofoot, onecolumn]{IEEEtran}
\usepackage{amsfonts}%
 \usepackage{setspace}
\usepackage{amssymb}%
\usepackage[dvips]{graphicx}%
\usepackage{amsmath}%
\usepackage[usenames]{color}
\usepackage{algorithmic}
\newtheorem{theorem}{Theorem}

\newtheorem{lemma}{Lemma}

\newtheorem{condition}{Condition}

\newcommand{\comment}[1]{}

\ifodd 0
\newcommand{\rev}[1]{{\color{blue}#1}} 
\newcommand{\com}[1]{\textbf{\color{red}(COMMENT: #1)}} 
\newcommand{\clar}[1]{\textbf{\color{green}(NEED more work: #1)}}

\else
\newcommand{\rev}[1]{#1}
\newcommand{\com}[1]{}
\newcommand{\clar}[1]{}
\fi

\begin{document}
\title{Online Learning of Rested and Restless Bandits} 
\author{\IEEEauthorblockN{Cem Tekin, Mingyan Liu\\
\thanks{Preliminary versions of this work appeared in Allerton 2010 and Infocom 2011.}} 
\IEEEauthorblockA{Department of Electrical Engineering and Computer Science\\
University of Michigan, Ann Arbor, Michigan, 48109-2122\\
Email: \{cmtkn, mingyan\}@umich.edu}
}

\maketitle

\begin{abstract}  
\rev{In this paper we study the online learning problem involving {\em rested} and {\em restless} multiarmed bandits
with multiple plays.  The system consists of a single player/user and a set of $K$ finite-state discrete-time Markov chains ({\em arms}) with unknown state spaces and statistics.  At each time step the player can play $M$, $M\leq K$, arms.
The objective of the user is to decide for each step which $M$ of the $K$ arms to play over a sequence of trials so as to maximize its long term reward.
The restless multiarmed bandit is particularly relevant to the application of opportunistic spectrum access (OSA), where a (secondary) user has access to a set of $K$ channels, each of time-varying condition as a result of random fading and/or certain primary users' activities.
} 
%

We first show that a logarithmic regret algorithm exists for the {\em rested} multiarmed bandit
problem.  We then construct an algorithm for the restless bandit problem which utilizes regenerative
cycles of a Markov chain and computes a sample mean based index
policy.  We show that under mild conditions on the state transition probabilities of the Markov chains this algorithm
achieves logarithmic regret uniformly over time, and that this regret bound is also optimal. 

\end{abstract}

\section{Introduction} \label{sec:intro}

In this paper we study the online learning problem involving {\em rested} and {\em restless} multiarmed bandits
with multiple plays.  The system consists of a single player/user and a set of $K$ finite-state discrete-time Markov chains (also referred to as {\em arms}) with unknown state spaces and statistics.  At each time step the player can play $M$, $M\leq K$, arms. Each arm played generates a reward depending on the state the arm is in when played.  
The state of an arm is only observed when it is played, and otherwise unknown to the user.  
The objective of the user is to decide for each step which $M$ of the $K$ arms to play over a sequence of trials so as to maximize its long term reward.  To do so it must use all its past actions and observations to essentially learn the quality of each arm (e.g., their expected rewards). 
\rev{We consider two cases, one with {\em rested} arms where the state of a Markov chain stays frozen unless it's played, the other with {\em restless} arms where the state of a Markov chain may continue to evolve (accordingly to a possibly different law) regardless of the player's actions.  } 

The above problem is motivated by the following opportunistic spectrum access (OSA) problem.  A (secondary)
user has access to a set of $K$ channels, each of time-varying condition as a result of random fading and/or certain primary
users' activities.  The condition of a channel is assumed to evolve as a Markov chain. 
At each time step, the secondary user (simply referred to as {\em the user} for the rest of the paper for there
is no ambiguity) senses or probes $M$ of the $K$ channels to find out their condition, and is allowed to use the channels in a way consistent with their conditions.  For instance, good channel conditions result in higher data rates or lower power for the user and so on.  In some cases channel conditions are simply characterized as being available and unavailable, and the user is allowed to use all channels sensed to be available.  
This is modeled as a reward collected by the user, the reward being a function of the state of the channel or the Markov chain.

\rev{The restless bandit model is particularly relevant to this application} because the state of each Markov chain evolves independently of the action of the user. 
The restless nature of the Markov chains follows naturally from the fact that channel conditions are governed by external factors like random fading, shadowing, and primary user activity. In the remainder of this paper a channel will also be referred to as an {\em arm}, the user as {\em player}, and probing a channel as {\em playing or selecting an arm}.

Within this context, the user's performance is typically measured by the notion of {\em regret}. It is defined as the difference between the expected reward that can be gained by an ``infeasible'' or ideal policy, i.e., a policy that requires either a priori knowledge of some or all statistics of the arms or hindsight information, and the expected reward of the user's policy.
The most commonly used infeasible policy is the {\em best single-action} policy, that is optimal among all policies that continue to play the same arm.  An ideal policy could play for instance the arm that has the highest expected reward (which requires statistical information but not hindsight).  This type of regret is sometimes also referred to as the {\em weak regret}, see e.g., work by Auer et al. \cite{auer2}.  In this paper we will only focus on this definition of regret.  Discussion on possibly stronger regret measures is given in Section \ref{sec:discussion}. 

This problem is a typical example of the tradeoff between {\em exploration} and {\em exploitation}.  On the one hand, the player needs to sufficiently explore all arms so as to discover with accuracy the set of best arms and avoid getting stuck playing an inferior one erroneously believed to be in the set of best arms.  On the other hand, the player needs to avoid spending too much time sampling the arms and collecting statistics and not playing the best arms often enough to get a high return.  


In most prior work on the class of multiarmed bandit problems, originally proposed by Robbins \cite{robbins2}, 
the rewards are assumed to be independently drawn from a fixed (but unknown) distribution.  
It's worth noting that with this iid assumption on the reward process, whether an arm is rested or restless is inconsequential for the following reasons.
Since the rewards are independently drawn each time, whether an unselected arm remains still or continues to change does not affect the reward the arm produces the next time it is played whenever that may be.
This is clearly not the case with Markovian rewards.
In the rested case, since the state is frozen when an arm is not played,
the state in which we next observe the arm is {\em independent} of how much time elapses before we play the arm again.   In the restless case, the state of an arm continues to evolve, 
thus the state in which we next observe it is now {\em dependent} on the amount of time that elapses between two plays of the same arm.  This makes the problem significantly more difficult.




Below we briefly summarize the most relevant results in the literature.
Lai and Robbins in  \cite{lai1} model rewards as single-parameter univariate densities 
and give a lower bound on the regret and construct policies that achieve this lower bound which are called {\em asymptotically efficient} policies. This result is extended by Anantharam et al. in \cite{anantharam2} to the case of playing more than one arm at a time. Using a similar approach Anantharam et al. in \cite{anantharam1} develops index policies that are asymptotically efficient for arms with rewards driven by finite, irreducible, aperiodic and rested Markov chains with identical state spaces and single-parameter families of stochastic transition matrices.  Agrawal in \cite{agrawal1} considers sample mean based index policies for the iid model that achieve $O(\log n)$ regret, where $n$ is the total number of plays.
Auer et al. in \cite{auer} also proposes sample mean based index policies for iid rewards with bounded support; these are derived from \cite{agrawal1}, but are simpler than those in \cite{agrawal1} and are not restricted to a specific family of distributions. These policies achieve logarithmic regret uniformly over time rather than asymptotically in time, but have bigger constant than that in \cite{lai1}.   
In \cite{tekin1} we showed that the index policy in \cite{auer} is order optimal for Markovian rewards drawn from rested arms but not restricted to single-parameter families, under some assumptions on the transition probabilities. 
\rev{Parallel to the work presented here, in \cite{hliu1} an algorithm was constructed that achieves logarithmic regret for the restless bandit problem.  The mechanism behind this algorithm however is quite different from what's presented here; this difference is discussed in more detail in Section \ref{sec:discussion}. } 

Other works such as \cite{liu2, anandkumar,krishnamachari} consider the iid reward case in a decentralized multiplayer setting; players selecting the same arms experience collision according to a certain collision model.  
We would like to mention another class of multiarmed bandit problems in which the statistics of the arms are known a priori and the state is observed perfectly; these are thus optimization problems rather than learning problems. The rested case is considered by Gittins \cite{gittins1} and the optimal policy is proved to be an index policy which at each time plays the arm with highest Gittins' index.  Whittle introduced the restless version of the bandit problem in \cite{whittle}.  The restless bandit problem does not have a known general solution though special cases may be solved.  For instance, a myopic policy is shown to be optimal when channels are identical and bursty in \cite{ahmad} for an OSA problem formulated as a restless bandit problem with each channel modeled as a two-state Markov chain (the Gilbert-Elliot model).

In this paper we first study the {rested} bandit problem with Markovian rewards. Specifically, we show that a straightforward extension of the UCB1 algorithm \cite{auer} to the multiple play case (UCB1 was originally designed for the case of a single play: $M=1$) results in logarithmic regret for restless bandits with Markovian rewards. 
%
We then use the key difference between rested and restless bandits to construct a regenerative cycle algorithm (RCA) that produces logarithmic regret for the restless bandit problem.  The construction of this algorithm allows us to use the proof of the rested problem as a natural stepping stone, and simplifies the presentation of the main conceptual idea. 

The work presented in this paper extends our previous results \cite{tekin1,tekin2} on single play to multiple plays ($M\geq 1$).  Note that this single player model with multiple plays at each time step is conceptually equivalent to the centralized (coordinated) learning by multiple players, each playing a single arm at each time step.  Indeed our proof takes this latter point of view for ease of exposition, and our results on logarithmic regret equally applies to both cases. 

The remainder of this paper is organized as follows. In Section
\ref{sec:problem} we present the problem formulation. In Section \ref{sec:rested-problem} we analyze a sample
mean based algorithm for the rested bandit problem.  In Section \ref{sec:restless-problem} we propose an algorithm
based on regenerative cycles that employs sample mean based indices and analyze its regret.
In Section \ref{sec:example} we numerically examine the performance of this algorithm in the case of an OSA problem with Gilbert-Elliot
channel model. In Section \ref{sec:discussion} we discuss possible improvements and compare our algorithm to other algorithms. 
Section \ref{sec:conc} concludes the paper.  


\section{Problem Formulation and Preliminaries} \label{sec:problem}

Consider $K$ arms (or channels) indexed by the set
$\mathcal{K}=\left\{1,2, \ldots, K\right\}$. The $i$th arm is
modeled as a discrete-time, irreducible and aperiodic Markov chain
with finite state space $S^i$. There is a stationary and positive
reward associated with each state of each arm. Let $r^i_x$ denote
the reward obtained from state $x$ of arm $i$, $x\in S^i$; this
reward is in general different for different states. \rev{Let
$P^i=\left\{p_{xy}^i, x,y \in S^i \right\}$ denote the transition
probability
matrix of the $i$-th arm, and 
$\boldsymbol{\pi}^i = \{\pi^i_x, x\in S^i \}$ the stationary
distribution of $P^i$.}

We assume the arms (the Markov chains) are mutually independent. In subsequent sections we will consider the rested and the restless cases separately. \rev{As mentioned in the introduction, the state of a rested arm changes according to $P^i$ only when it is played and remains frozen otherwise.  By contrast, the state of a restless arm changes according to $P^i$ regardless of the user's actions.} All the assumptions in this section applies to both types of arms. We note that the rested model is a special case of the restless model, but our development under the restless model follows the rested model\footnote{\rev{
In general a restless arm may be given by two transition probability matrices, an active one ($P^i$) and a passive one ($Q^i$). 
The first describes the state evolution when it is played and the second the state evolution when it is not played. 
When an arm models channel variation, $P^i$ and $Q^i$ are in general assumed to be
the same as the channel variation is uncontrolled.  In the context of online learning we shall see 
that the selection of $Q^i$ is irrelevant; indeed the arm does not even have to be Markovian when
it's in the passive mode.  More is discussed in Section \ref{sec:discussion}.}}. 

Let $(P^i)'$ denote the {\em adjoint} of $P^i$ on $l_2(\pi)$ where
\begin{eqnarray*}
(p^i)'_{xy}=(\pi^i_y p^i_{yx})/\pi^i_x, ~ \forall x,y \in S^i,
\end{eqnarray*}
and $\hat{P}^i=(P^i)'P$ denotes the {\em multiplicative symmetrization} of $P^i$. We will assume that the $P^i$s are such that $\hat{P}^i$s are irreducible.  \rev{To give a sense of how weak or strong this assumption is, we first note that this is a weaker condition than assuming the Markov chains to be reversible.
In addition, we note that one condition that guarantees the $\hat{P}^i$s are irreducible is $p_{xx}>0, \forall x\in S^i, \forall i$. This assumption thus holds naturally for our main motivating application, as it's possible for channel condition to remain the same over a single time step (especially if the unit is sufficiently small).  It also holds for a very large class of Markov chains and applications in general.  Consider for instance a queueing system scenario where an arm denotes a server and the Markov chain models its queue length, in which it is possible for the queue length to remain the same over one time unit. 
} 


The mean reward of arm $i$, denoted by $\mu^i$, is the expected reward of arm $i$ under its stationary distribution: 
\begin{eqnarray}
\mu^i=\displaystyle\sum_{x \in S^i} r^i_x \pi^i_x ~.
\end{eqnarray}
%

\rev{Consistent with the discrete time Markov chain model, we will assume that the player's actions occur in discrete time steps. Time is indexed by $t$, $t=1, 2, \cdots$.  We will also frequently refer to the time interval $(t-1, t]$ as time slot $t$.} 
The player plays $M$ of the $K$ arms at each time step. 

\rev{Throughout the analysis we will make the additional assumption that the mean reward of arm $M$ is strictly greater than the mean reward of arm $M+1$, i.e., we have $\mu^1 \geq \mu^2 \geq \cdots \geq \mu^M > \mu^{M+1} \geq \cdots \geq \mu^K$. For rested arms this assumption simplifies the presentation and is not necessary, i.e., results will hold for $\mu^M \geq \mu^{M+1}$. However, for restless arms the strict inequality between $\mu^M$ and $\mu^{M+1}$ is needed because otherwise there can be a large number of arm switchings between the $M$-th and the $(M+1)$-th arms (possibly more than logarithmic).  Strict inequality will prevent this from happening. 
We note that this assumption is not in general restrictive; in our motivating application distinct channel conditions typically means different data rates. Possible relaxation of this condition is given in Section \ref{sec:discussion}.} 

We will refer to the set of arms $\left\{1, 2, \cdots, M\right\}$ as the $M$-best arms and say that each arm in this set is {\em optimal} while referring to the set $\left\{M+1, M+2, \cdots, K\right\}$ as the $M$-worst arms and say that each arm in this set is {\em suboptimal}.  

For a policy $\alpha$ we define its regret $R^\alpha(n)$ as the difference between the expected total reward that can be obtained by only playing the $M$-best arms and the expected total reward obtained by policy $\alpha$ up to time $n$.   
Let \rev{$A^{\alpha}(t)$} denote the set of arms selected by policy $\alpha$ at $t$, $t=1, 2, \cdots$, and $x_{\alpha}(t)$ the state of arm $\alpha(t) \in A^{\alpha}(t)$ at time $t$.  Then we have
\begin{eqnarray}
R^\alpha(n)=n \sum_{j=1}^M \mu^{j}-E^\alpha\left[\displaystyle\sum_{t=1}^n \sum_{\alpha(t) \in A^{\alpha}(t)} r^{\alpha(t)}_{x_{\alpha}(t)} \right] ~.
\end{eqnarray}
The objective is to examine how the regret $R^\alpha(n)$ behaves as a function of $n$ for a given policy $\alpha$ and to construct a policy whose regret is order-optimal, through appropriate bounding.  As we will show and as is commonly done, the key to bounding $R^\alpha(n)$ is to bound the expected number of plays of any suboptimal arm.
%

Our analysis utilizes the following known results on Markov
chains; the proofs are not reproduced here for brevity.  The
\rev{first result is due to Lezaud} \cite{lezaud} that bounds the
probability of a large deviation from the stationary distribution.
\begin{lemma}\label{lemma3}
[Theorem 3.3 from \cite{lezaud}] Consider a finite-state, irreducible Markov chain $\left\{X_t\right\}_{t \geq 1}$ with state space $S$, matrix of transition probabilities $P$, an initial distribution $\mathbf{q}$ and stationary distribution $\mathbf{\pi}$. Let $N_{\mathbf{q}}=\left\|(\frac{q_x}{\pi_x}, x\in S)\right\|_2$. Let $\hat{P}=P'P$ be the multiplicative symmetrization of $P$ where $P'$ is the adjoint of $P$ on $l_2(\pi)$. Let $\epsilon=1-\lambda_2$, where $\lambda_2$ is the second largest eigenvalue of the matrix $\hat{P}$. $\epsilon$ will be referred to as the eigenvalue gap of $\hat{P}$. Let $f:S\rightarrow \mathbb{R}$ be such that $\sum_{y \in S} \pi_y f(y) =0$, $\left\|f\right\|_{\infty} \leq 1$ and $0<\left\|f\right\|^2_2 \leq 1$. If $\hat{P}$ is irreducible, then for any positive integer $n$ and all $0 < \gamma \leq 1$
\begin{eqnarray*}
P\left(\frac{\sum_{t=1}^n f(X_t)}{n} \geq \gamma \right) \leq N_q \exp\left[-\frac{n \gamma^2 \epsilon}{28}\right] ~.
\end{eqnarray*}
\end{lemma}

The second result is due to Anantharam et al., which can be found in \cite{anantharam1}. 
\begin{lemma}\label{lemma:anantharam}
[Lemma 2.1 from \cite{anantharam1}] Let $Y$ be an irreducible aperiodic Markov chain with a state space $S$, transition probability matrix $P$,  an initial distribution that is non-zero in all states, and 
a stationary distribution $\{\pi_x\}, \forall x\in S$. Let $F_t$ be the $\sigma$-field generated by random variables $X_1,X_2,...,X_t$ where $X_t$ corresponds to the state of the chain at time $t$.  Let $G$ be a $\sigma$-field independent of $F=\vee_{t \geq 1} F_t$, the smallest $\sigma$-field containing $F_1, F_2, ...$. Let $\tau$ be a stopping time with respect to the increasing family of $\sigma$-fields $\left\{G \vee F_t, t \geq 1\right\}$. Define  $N(x,\tau)$ such that
\begin{eqnarray*}
N(x,\tau)=\displaystyle\sum_{t=1}^\tau I(X_t = x).
\end{eqnarray*}
Then $\forall \tau$ such that $E\left[\tau\right]<\infty$, we have
\begin{eqnarray}
\left|E\left[N(x,\tau)\right]-\pi_xE\left[\tau\right]\right|\leq C_{P},
\end{eqnarray}
where $C_{P}$ is a constant that depends on $P$.
\end{lemma}

The third result is due to Bremaud, which can be found in \cite{bremaud}.
\begin{lemma}
If $\left\{X_n\right\}_{n \geq 0}$ is a positive recurrent
homogeneous Markov chain with state space $S$, stationary
distribution $\pi$ and $\tau$ is a stopping time that is finite
almost surely for which $X_{\tau}=x$ then for all $y \in S$
\begin{eqnarray*}
E\left[ \sum_{t=0}^{\tau-1} I(X_t=y) | X_0=x \right] = E[\tau | X_0=x]\pi_y ~.
\end{eqnarray*}
\label{lastlemma}
\end{lemma}

\rev{The following notations are frequently used
throughout the paper: $\beta= \sum_{t=1}^{\infty} 1/t^{2}$, $\pi^i_{\min} = \min_{x \in S^i} \pi^i_x$,
$\pi_{\min} = \min_{i \in \mathcal{K}} \pi^i_{\min}$,
$r_{\max}=\max_{x \in S^i, i \in \mathcal{K}} r^i_x$,
$S_{\max}=\max_{i \in \mathcal{K}} |S^i|$, $\hat{\pi}_{\max} =
\max_{x \in S^i, i \in \mathcal{K}} \left\{\pi^i_x,
1-\pi^i_x\right\}$, $\epsilon_{\min}= \min_{i \in \mathcal{K}}
\epsilon^i$, where $\epsilon^i$ is the eigenvalue gap ({the difference between 1 and 
the second largest eigenvalue}) 
of the multiplicative symmetrization of the transition probability matrix
of the $i$th arm, 
and $\Omega^i_{\max} = \max_{x, y \in S^i}
\Omega^i_{x, y}$, where $\Omega^i_{x, y}$ is the mean hitting time
of state $y$ given the initial state $x$ for arm $i$ for $P^i$.}

In the next two sections we present algorithms for the {rested}
and {restless} problems, referred to as the {\em upper
confidence bound - multiple plays} (UCB-M) and the {\em regenerative cycle
algorithm - multiple plays} (RCA-M), respectively, and analyze their
regret.

\section{\rev{Analysis of} the Rested Bandit Problem with Multiple Plays} \label{sec:rested-problem}

In this section we \rev{show that there exists an algorithm that achieves logarithmic regret uniformly over time for the 
rested bandit problem with Markovian reward and multiple plays.}  We present such an algorithm, called {\em the upper
confidence bound - multiple plays} (UCB-M), which is a straightforward extension of UCB1 from \cite{auer}.  
This algorithm plays $M$ of the $K$ arms with the highest indices with a modified
exploration constant $L$ instead of $2$ in \cite{auer}. Throughout
our discussion, we will consider a horizon of $n$ time slots.
\rev{For simplicity of presentation we will view a single player
playing multiple arms at each time as multiple coordinated players each
playing a single arm at each time.} In other
words we consider $M$ players indexed by $1,2,\cdots,M$, each
playing a single arm at a time.  Since in this case information
is centralized, collision is completely avoided among the players, i.e., at each time
step an arm will be played by at most one player.

Below we summarize a list of notations used in this section.

\begin{itemize}
\item $A(t)$: the set of arms played at time $t$ \rev{(or in slot $t$)}.
\item $T^i(t)$: total number of times (slots) arm $i$ is played \rev{up to} the end of slot $t$.
\item $T^{i,j}(t)$: total number of times (slots) player $j$ played arm $i$ \rev{up to} the end of slot $t$.
\item $\bar{r}^i (T^i(t))$: sample mean of the rewards observed from the first $T^i(t)$ plays \rev{of} arm $i$.
\end{itemize}

As shown in Figure \ref{fig:UCBM}, UCB-M selects $M$ channels with
the highest indices at each time step and updates the indices
according to the rewards observed. The index given on line 4 of
Figure \ref{fig:UCBM} depends on the sample mean reward and an
exploration term which reflects the relative uncertainty about the
sample mean of an arm. We call $L$ in the exploration term
{\em the exploration constant}. The exploration term grows
logarithmically when the arm is not played in order to guarantee
that sufficient samples are taken from each arm to approximate the
mean reward.

\begin{figure}[htb]
\fbox {
\begin{minipage}{\columnwidth}
{The Upper Confidence Bound - Multiple Plays (UCB-M): }
\begin{algorithmic}[1]
\STATE {Initialize: Play each arm $M$ times in the first $K$ slots}
\WHILE {$t\geq K$}
\STATE{$\bar{r}^i (T^i(t)) =\frac{r^i(1)+r^i(2)+...+r^i(T^i(t))}{T^i(t)}, ~\forall i$}
\STATE{calculate index: $g^i_{t,T^i(t)}= \bar{r}^i(T^i(t)) + \sqrt{\frac{L\ln t}{T^i(t)}}, ~\forall i$}
\STATE{$t:=t+1$}
\STATE{play $M$ arms with the highest indices, update $r^j(t)$ and $T^j(t)$.}
\ENDWHILE
\end{algorithmic}
\end{minipage}
} \caption{pseudocode for the UCB-M algorithm.} \label{fig:UCBM}
\end{figure}

\rev{To upper bound the regret of the above algorithm logarithmically, we proceed as follows.  We begin by relating the
regret to the expected number of plays of the arms and then show
that each suboptimal arm is played at most logarithmically in
expectation.  These steps are illustrated in the following lemmas. 
Most of these lemmas are established under the following condition on the arms. 
\begin{condition}\label{cond:1}
All arms are finite-state, irreducible, aperiodic Markov chains whose transition probability matrices have
irreducible multiplicative symmetrizations and $r^i_x >0$,
$\forall i \in \mathcal{K}$, $\forall x \in S^i$. 
\end{condition}
} 


\begin{lemma} \label{lemma:rested1}
Assume that all arms are finite-state, irreducible, aperiodic, rested Markov chains. Then \rev{using UCB-M we have:} 
\begin{eqnarray}
\left|R(n)-\left(n\sum_{j=1}^M \mu^{j}- \sum_{i=1}^K \mu^i E[T^i(n)]\right)\right| \leq C_{\mathbf{S,P,r}},
\end{eqnarray}
where $C_{\mathbf{S,P,r}}$ is a constant that depends on the state \rev{spaces}, rewards, \rev{and} transition probabilities but not on time.
\end{lemma}
\begin{proof}
see Appendix \ref{app:B}.
\end{proof}

\begin{lemma} \label{lemma:rev1}
\rev{Assume Condition \ref{cond:1} holds and all arms are rested.}  
\rev{Under UCB-M} with $L \geq 112 S^2_{\max} r^2_{\max}
\hat{\pi}^2_{\max} /\epsilon_{\min}$, for any suboptimal arm $i$, \rev{we have} 
\begin{eqnarray}
E[T^i(n)] \leq M + \frac{4L\ln n}{(\mu^{M}-\mu^i)^2} + \sum_{j=1}^M \frac{(|S^i|+|S^{j}|)\beta}{\pi_{\min}} \nonumber 
\end{eqnarray}
\end{lemma}
\begin{proof}
see Appendix \ref{app:C}.
\end{proof}

\begin{theorem} \label{thm:mainrested}
\rev{Assume Condition \ref{cond:1} holds and all arms are rested.}  
\rev{With} constant $L \geq 112 S^2_{\max} r^2_{\max}
\hat{\pi}^2_{\max} /\epsilon_{\min}$ the regret \rev{of UCB-M} is upper bounded
by
\begin{eqnarray}
R(n) \leq 4L\ln n \sum_{i>M} \frac{(\mu^{1}-\mu^i)}{(\mu^{M}-\mu^i)^2} + \sum_{i>M} (\mu^{1}-\mu^i) \left(M + \sum_{j=1}^M C_{i,j} \right) + C_{\mathbf{S,P,r}} , \nonumber \\
\end{eqnarray}
where
$C_{i,j} =\frac{(|S^i|+|S^{j}|)\beta} {\pi_{\min}}$.
\end{theorem}
\begin{proof}
\begin{eqnarray}
n \sum_{j=1}^M \mu^{j} - \sum_{i=1}^K \mu^i E[T^i(n)] &=& \sum_{j=1}^M \sum_{i=1}^K \mu^{j} E[T^{i,j}(n)] - \sum_{j=1}^M \sum_{i=1}^K \mu^i E[T^{i,j}(n)] \nonumber \\
&=& \sum_{j=1}^M \sum_{i>M} (\mu^{j}- \mu^i) E[T^{i,j}(n)] \leq \sum_{i>M} (\mu^{1}-\mu^i) E[T^i(n)]. \nonumber 
\end{eqnarray}
Thus,
\begin{eqnarray}
R(n) &\leq& n \sum_{j=1}^M \mu^{j} - \sum_{i=1}^K \mu^i E[T^i(n)] + C_{\mathbf{S,P,r}} \label{eqn:ver11} \\
&\leq& \sum_{i>M} (\mu^{1}-\mu^i) E[T^i(n)] + C_{\mathbf{S,P,r}} \nonumber \\
&\leq& \sum_{i>M} (\mu^{1}-\mu^i) \left(M + \frac{4L\ln n}{(\mu^{M}-\mu^i)^2} + \sum_{j=1}^M \frac{(|S^i|+|S^{j}|)\beta}{\pi_{\min}}\right) + C_{\mathbf{S,P,r}} \label{eqn:ver12} \\
&=& 4L\ln n \sum_{i>M} \frac{(\mu^{1}-\mu^i)}{(\mu^{M}-\mu^i)^2} + \sum_{i>M} (\mu^{1}-\mu^i) \left(M + \sum_{j=1}^M C_{i,j} \right) + C_{\mathbf{S,P,r}}, \nonumber
\end{eqnarray}
where (\ref{eqn:ver11}) follows from Lemma \ref{lemma:rested1} and  (\ref{eqn:ver12}) follows from Lemma \ref{lemma:rev1}.
\end{proof}

\rev{The above theorem says that provided that $L$ satisfies the stated sufficient condition, UCB-M results in logarithmic regret for the rested problem.  This sufficient condition does require certain knowledge on the underlying Markov chains.  This requirement may be removed if the value of $L$ is adapted over time.  More is discussed in Section \ref{sec:discussion}.} 

\section{\rev{Analysis of} the Restless Bandit Problem with Multiple Plays} \label{sec:restless-problem}

\rev{In this section we study the restless bandit problem.  We construct an algorithm called the} {\em regenerative cycle algorithm - multiple plays} (RCA-M), and prove that this algorithm guarantees logarithmic regret uniformly over time under the \rev{same mild assumptions on the state transition probabilities as in the rested case}. 
RCA-M is a multiple plays extension of RCA first introduced in \cite{tekin2}.  \rev{Below we first present the key conceptual idea behind RCA-M, followed by a more detailed pseudocode.  We then prove the logarithmic regret result. } 

As the name suggests, RCA-M operates in regenerative cycles.  
In essence RCA-M uses the observations from sample paths
\rev{within} regenerative cycles to estimate the \rev{sample
mean of an arm in the form of an index similar to that used in UCB-M} while discarding the rest of the
observations \rev{(only for the computation of the index, but they are added to the total reward)}.  
\rev{Note that the
rewards from the discarded observations are collected but are not
used to make decisions.} The reason behind such a construction has
to do with the restless nature of the arms. Since each arm
continues to evolve according to the Markov chain regardless of
the user's action, the probability distribution of the reward we
get by playing an arm is a function of the amount of time that has
elapsed since the last time we played the same arm.  Since the
arms are not played continuously, the sequence of observations
from an arm which is not played consecutively does not correspond
to a discrete time homogeneous Markov chain. While this certainly
does not affect our ability to collect rewards, it becomes hard to
analyze the estimated quality (the index) of an arm calculated
based on rewards collected this way.

However, if instead of the actual sample path of observations from an arm, we limit ourselves to a sample path constructed (or rather stitched together) using only the observations from regenerative cycles, then this sample path essentially has the same statistics as the original Markov chain due to the renewal property and one can now use the sample mean of the rewards from the regenerative sample paths to approximate the mean reward under stationary distribution.  

\rev{Under RCA-M each player maintains a block structure; a block consists of a certain number of slots}. 
\rev{Recall that as mentioned earlier, even though our basic model is one of single-player multiple-play, our description is in the equivalent form of multiple coordinated players each with a single play.}  Within a block a player plays the same arm continuously till a certain pre-specified state (say $\gamma^i$) is observed. Upon this observation \rev{the arm} enters a regenerative cycle and the player continues to play the same arm till state $\gamma^i$ is observed for the second time, which denotes the end of the block. Since $M$ arms are played (by $M$ players) simultaneously in each slot, different blocks overlap in time. Multiple blocks may or may not start or end at the same time. \rev{In our analysis below blocks will be ordered; they are ordered according to their start time. If multiple blocks start at the same time then the ordering among them is randomly chosen.} 

\begin{figure}[ht]
\begin{center}
\includegraphics[width=6in]{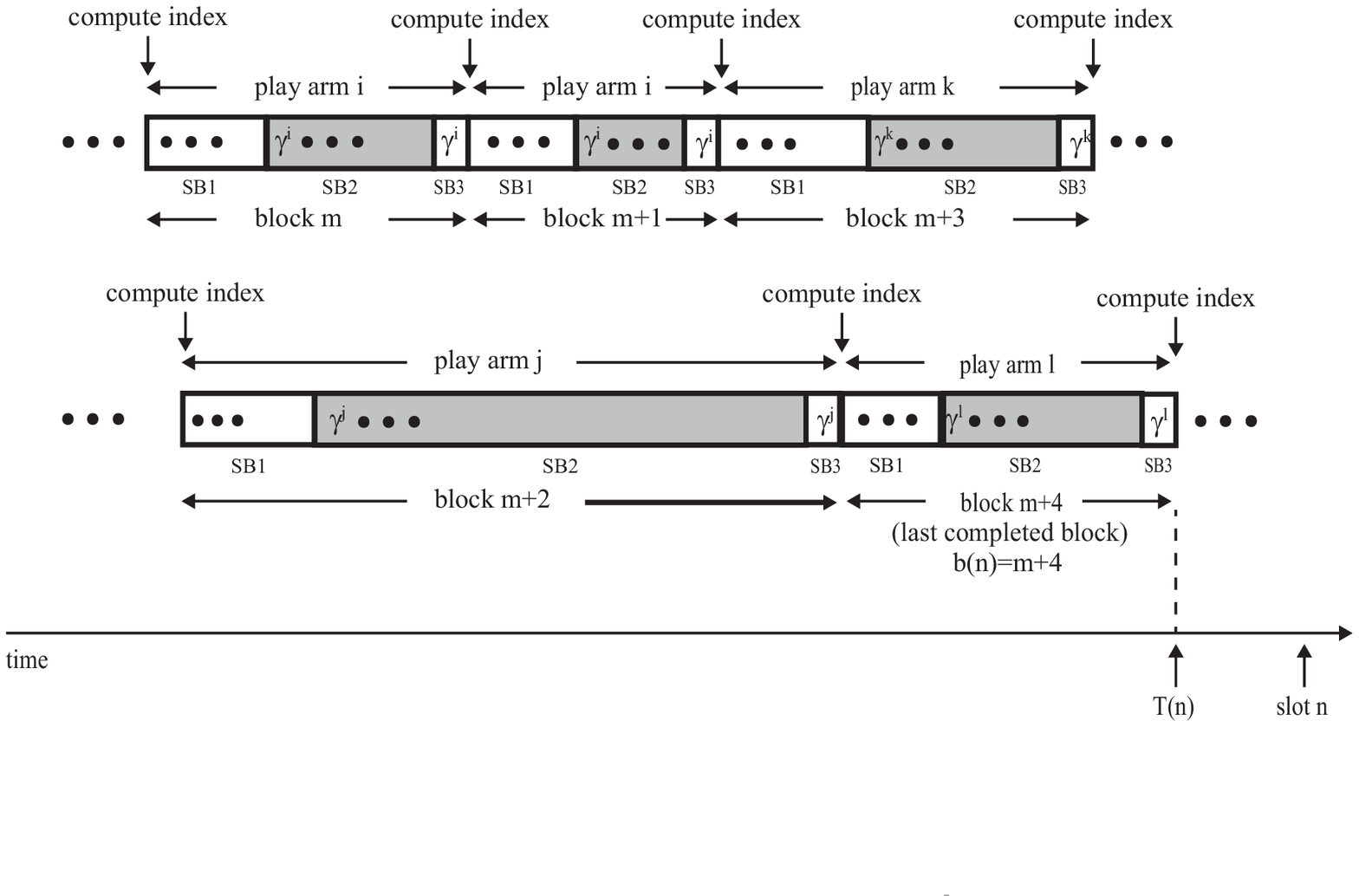}
\vspace{-0.8in}
\caption{Example realization of RCA-M with $M=2$ for a period of $n$ slots}
\label{fig:RCAM2}
\end{center}
\end{figure}

For the purpose of index computation and subsequent analysis, each block is further broken into three sub-blocks (SBs).  SB1 consists of all time slots from the beginning of the block to right before the first visit to $\gamma^i$; SB2 includes all time slots from the first visit to $\gamma^i$ up to but excluding the second visit to state $\gamma^i$; SB3 consists of a single time slot with the second visit to $\gamma^i$. 
\rev{Figure \ref{fig:RCAM2} shows an example sample path of the operation of RCA-M.  The block structure of two players are shown in this example; the ordering of the blocks is also shown. } 


The key to the RCA-M algorithm is for each arm to single out only observations within SB2's in each block and virtually assemble them. Throughout our discussion, we will consider a horizon of $n$ time slots.  A list of notations used is summarized as follows:
\begin{itemize}
\item $A(t)$: the set of arms played at time $t$ \rev{(or in time slot $t$)}.
\item $\gamma^i$: \rev{the} state that \rev{determines} the regenerative cycles for arm $i$.
\item $\alpha(b)$: the arm played in \rev{the $b$-th} block.
\item $b(n)$: \rev{the} total number of completed blocks \rev{by all players} up to time $n$.
\item $T(n)$: \rev{the} time at the end of the last completed block across all arms \rev{(see Figure \ref{fig:RCAM2})}.
\item $T^i(n)$: \rev{the} total number of times (slots) arm $i$ is played up to the last completed block of arm $i$ up to time $T(n)$.
\item $T^{i,j}(n)$: \rev{the} total number of times (slots) arm $i$ is played by user $j$ up to the last completed block of arm $i$ up to time $T(n)$
\item $B^i(b)$: \rev{the} total number of blocks within the first completed $b$ blocks in which arm $i$ is played.

\item $X^i_1(b)$: \rev{the} vector of observed states from SB1 of the \rev{$b$-th} block in which arm $i$ is played; 
\rev{this vector} is empty if the first observed state is $\gamma^i$.
\item $X^i_2(b)$: \rev{the} vector of observed states from SB2 of the \rev{$b$-th} block in which arm $i$ is played;
\item $X^i(b)$: \rev{the} vector of observed states from the \rev{$b$-th} block in which arm $i$ is played. Thus we have
$X^i(b)=[X^i_1(b), X^i_2(b), \gamma^i] $.

\item $t(b)$: time at the end of block $b$;

\item $T^i(t(b))$: \rev{the} total number of time slots arm $i$ is played up to the last completed block of arm $i$ \rev{within} time $t(b)$.

\item $t_2(b)$: \rev{the} total number of time slots that \rev{lie within} at least one SB2 in a completed block of any arm up to \rev{and including} block $b$.


\item $r^i(t)$: the reward from arm $i$ \rev{upon its $t$-th play}, counting only those \rev{plays} during an SB2.

\item $T^i_2(t_2(b))$: \rev{the} total number of time slots arm $i$ is played during SB2's up to \rev{and including} block $b$.

\item $O(b)$: \rev{the set of arms that are {\em free} to be selected by some player $i$ upon its completion of the $b$-th block; these are arms that are currently not being played by other players (during time slot $t(b)$), and the arms whose blocks are completed at time $t(b)$.}  
\end{itemize}


RCA-M computes and updates the value of an {\em index} $g^i$ for each arm $i$ in \rev{the set $O(b)$} at the end of block $b$ based on the total reward obtained from arm $i$ during all \rev{SB2's} as follows:
\begin{eqnarray}
g^i_{t_2(b), T^i_2(t_2(b))}= \bar{r}^i(T^i_2(t_2(b))) + \sqrt{\frac{L \ln t_2(b)}{T^i_2(t_2(b))}}, \label{eqn1redifine}
\end{eqnarray}
where $L$ is a constant, and
\begin{eqnarray}
\rev{\bar{r}^i(T^i_2(t_2(b)))}=\frac{r^i(1)+r^i(2)+...+r^i(T^i_2(t_2(b)))}{T^i_2(t_2(b))} \nonumber
\end{eqnarray}
denotes the sample mean of the reward collected during SB2. 
\rev{Note that this is the same way the index is computed under UCB-M if we only consider SB2's.} 
\rev{Its also worth noting that under RCA-M rewards are also collected during SB1's and SB3's}.  However, the computation of the indices only relies on SB2. 
The pseudocode of RCA-M is given in Figure \ref{fig:RCAM}. 

\rev{Due to the regenerative nature of the Markov chains,} the rewards used in \rev{the} computation of the index of an arm can be viewed as rewards from a rested arm with the same transition matrix as the active transition matrix of the restless arm. 
However, \rev{to prove the} existence of a logarithmic upper bound \rev{on the regret} for restless arms \rev{remains a} non-trivial \rev{task} since the blocks may be arbitrarily long and the frequency of arm selection depends on the length of the blocks. 

\rev{In the analysis that follows,} we first show that the expected number of blocks in which a suboptimal arm is played is at most logarithmic by applying the result in Lemma \ref{lemma-key} \rev{that compares} the indices of arms in slots where an arm is selected. Using this result we \rev{then} show that the expected number of blocks in which a suboptimal arm is played is at most logarithmic in time. Using irreducibility of the arms the expected block length is finite, thus the number of time slots in which a suboptimal arm is played is finite. Finally, we show that the regret due to arm switching is at most logarithmic.


%

\begin{figure}[ht]
\fbox {
\begin{minipage}{\columnwidth }
{The Regenerative Cycle Algorithm - Multiple Plays (RCA-M):}\hspace{-30pt}
\begin{algorithmic}[1]
\STATE {Initialize: $b=1, t=0, t _2= 0, T^i_2=0, r^i=0, I^i_{SB2}=0, I^i_{IN}=1, \forall i=1,\cdots,K$, ${A}= \emptyset$}
\STATE{\rev{~~~//$I^i_{IN}$ indicates whether arm $i$ has been played at least once}} 
\STATE{\rev{~~~//$I^i_{SB2}$ indicates whether arm $i$ is in an SB2 sub-block}} 
\WHILE{(1)}
\FOR{$i=1$ to $K$}
\IF{$I^i_{IN}=1$ and $|{A}|<M$} 
\STATE{${A} \leftarrow {A} \cup \left\{i\right\}$ \rev{~~ //arms never played is given priority to ensure all arms are sampled initially} }
\ENDIF
\ENDFOR
\rev{
\IF{$|{A}|<M$}
\STATE{Add to ${A}$ the set $\left\{i: 
g^i \mbox{ is one of the } M-|{A}| \mbox{ largest among } \{g^k, k \in\{1, \cdots, K\}-{A}\} \right\}$} 
\STATE{~~~//for arms that have been played at least once, those with the largest indices are selected} 
\ENDIF
} 
\FOR{$i \in {A}$}
\STATE{play arm $i$; denote state observed by $x^i$}
\IF{$I^i_{IN}=1$}
\STATE{$\gamma^i=x^i$, $T^i_2:=T^i_2+1$, $r^i:=r^i+r^i_{x^i}$, $I^i_{IN}=0$, $I^i_{SB2}=1$}
\STATE{\rev{~~~//the first observed state becomes the regenerative state; the arm enters SB2}} 
\ELSIF{$x^i \neq \gamma^i$ and $I^i_{SB2}=1$}
\STATE{$T^i_2:=T^i_2+1$, $r^i:=r^i+r^i_{x^i}$}
\ELSIF{$x^i = \gamma^i$ and $I^i_{SB2}=0$}
\STATE{$T^i_2:=T^i_2+1$, $r^i:=r^i+r^i_{x^i}$, $I^i_{SB2}=1$}
\ELSIF{$x^i = \gamma^i$ and $I^i_{SB2}=1$}
\STATE{\rev{$r^i:=r^i+r^i_{x^i}$,} $I^i_{SB2}=0$, ${A} \leftarrow {A} - \left\{i\right\}$}
\ENDIF
\ENDFOR
\STATE{$t := t+1$, $t_2 := t_2+\min{\left\{1, \sum_{i\in S} I^i_{SB2}\right\}}$
\rev{~~//$t_2$ is only accumulated if at least one arm is in SB2}}
\FOR{$i=1$ to $K$}
\STATE{$g^i=\frac{r^i}{T^i_2}+\sqrt{\frac{L \ln t_2}{T^i_2}}$}
\ENDFOR
\ENDWHILE
\end{algorithmic}
\end{minipage}
}
\caption{Pseudocode of RCA-M} \label{fig:RCAM}
\end{figure}

\comment{
\hspace{0.2in}\begin{framed}
\noindent\textbf{Regenerative Cycle Algorithm (RCA)} \\

\noindent\textbf{Initialization:} $b=1$, $t_2=0$, $T^i_2(t_2)=0, \forall i=1,\cdots,K$\\

for ($b \leq K$) \\
\textbf{Initialization:} Play arm $b$. Set $\gamma^b$ to the first state observed. \\
    \indent $t_2$=$t_2+1$, $T^b_2(t_2)=T^b_2(t_2)+1$.\\
\textbf{Loop:}\\
    \textbf{if} observed state is not $\gamma^b$ keep playing the same arm.\\
    \indent  $t_2$=$t_2+1$, $T^b_2(t_2)=T^b_2(t_2)+1$\\
    \textbf{else if} observed state is $\gamma^b$  \\
    \indent  $b=b+1$\\

while ($b>K$) \\
In block $b$ play arm with the highest index according to the indices calculated at the end of block $b-1$.\\
\textbf{Loop:}\\
    \textbf{if} $\gamma^{\alpha(b)}$ is not observed up to current time in block $b$\\
    \indent wait\\
    \textbf{else if} $\gamma^{\alpha(b)}$ is observed once in block $b$ up to current time\\
    \indent $t_2$=$t_2+1$, $T^{\alpha(b)}_2(t_2)=T^{\alpha(b)}_2(t_2)+1$\\
    \textbf{else if} $\gamma^{\alpha(b)}$ is observed for the second time in block $b$\\
    \indent update the indices,$b=b+1$\\
\end{framed}
}


We bound the expected number of plays from a suboptimal arm.
\begin{lemma} \label{lemma:restless1}
\rev{Assume Condition \ref{cond:1} holds and all arms are restless.} 
\rev{Under} RCA-M with a constant $L \geq 112 S^2_{\max} r^2_{\max} \hat{\pi}^2_{\max} /\epsilon_{\min}$, we have
\begin{eqnarray}
&& \sum_{i>M} (\mu^{1}-\mu^i) E[T^i(n)] \leq  4L \sum_{i>M} \frac{(\mu^{1}-\mu^i) D_i \ln n}{(\mu^{M}-\mu^i)^2} +  \sum_{i>M} (\mu^{1}-\mu^i) D_i \left(1 + M \sum_{j=1}^M   C_{i,j}\right) ~, \nonumber
\end{eqnarray}
where
\begin{eqnarray*}
C_{i,j} &=& \frac{(|S^i|+|S^{j}|)\beta}{\pi_{\min}}, ~~ \beta = \sum_{t=1}^{\infty} t^{-2}, ~~ D_i = \left(\frac{1}{\pi^i_{\min}} + M^i_{\max} +1\right).
\end{eqnarray*}
\end{lemma}
\begin{proof}
see Appendix \ref{app:E}.
\end{proof}

We now state the main result of this section.
\begin{theorem} \label{thm:restless1}
\rev{Assume Condition \ref{cond:1} holds and all arms are restless.} 
\rev{With} constant $L \geq 112 S^2_{\max} r^2_{\max} \hat{\pi}^2_{\max} /\epsilon_{\min}$ the regret 
\rev{of RCA-M} is upper bounded by
\begin{eqnarray*}
R(n) &<& 4L \ln n \sum_{i>M} \frac{1}{(\mu^{M}-\mu^i)^2}\left((\mu^{1}-\mu^i) D_i + E_i\right) \\
&+& \sum_{i>M} \left((\mu^{1}-\mu^i)D_i+E_i\right) \left(1+ M \sum_{j=1}^M C_{i,j}\right) +F
\end{eqnarray*}
where
\begin{eqnarray*}
C_{i,j} &=& \frac{(|S^i|+|S^{j}|)\beta} {\pi_{\min}}, ~~ \beta = \sum_{t=1}^{\infty} t^{-2}\\
D_i &=& \left(\frac{1}{\pi^i_{\min}} + M^i_{\max} +1\right),\\
E_i &=& \mu^i(1+M^i_{\max}) + \sum_{j=1}^M \mu^{j}M^{j}_{\max},\\
F &=& \sum_{j=1}^M \mu^{j}\left(\frac{1}{\pi_{\min}} + \max_{i \in \mathcal{K}} M^i_{\max} +1\right).
\end{eqnarray*}
\end{theorem}
\begin{proof}
see Appendix \ref{app:F}. 
\end{proof}

\rev{Theorem \ref{thm:restless1} suggests} that given minimal information about the arms such as an
upper bound for $S^2_{\max} r^2_{\max} \hat{\pi}^2_{\max}
/\epsilon_{\min}$ the player can guarantee logarithmic regret by
choosing an $L$ in RCA-M that satisfies the \rev{stated} condition. 
\rev{As the rested case, this requirement on $L$ can be completely removed if the value of $L$ is adapted over time; more is discussed in Section \ref{sec:discussion}. } 

\comment{\rev{We say that the logarithmic bound in $n$ is {\em order
optimal} in the worst-case sense. Note that even if the transition matrices
are known a priori, the optimal policy for the restless bandits is not known in general
(solutions in some special cases are sometimes available).  
Therefore we compare the performance of our
algorithm with the worst-case optimal, i.e., minimum of the reward
of the optimal policy over all environments with the same active
transition matrices and arbitrary passive transition matrices. \com{Let's talk about this last sentence...}}
\com{Also should we talk about regret w.r.t. optimal policy when it's available?} }

\comment{i.e., no better bound than $\ln n$ is possible (however a
better constant is possible). This follows from the fact that the
rested bandit problem is a special case of the restless problem
and in \cite{anantharam1} it is shown that the best order is
logarithmic for the rested problem.}

We conjecture that the order optimality of RCA-M holds when it is
used with any index policy that is order optimal for the rested
bandit problem. Because of the use of regenerative cycles in
RCA-M, the observations used to calculate the indices can be in
effect treated as coming from rested arms. Thus an approach
similar to the one \rev{used} in the proof of Theorem \ref{thm:restless1} can
be used to prove order optimality \rev{of combinations of RCA-M and other index policies}.

\section{An Example for OSA: Gilbert-Elliot Channel Model} \label{sec:example}

In this section we simulate RCA-M under the commonly used
Gilbert-Elliot channel model where each channel has two states,
{\em good} and {\em bad} (or 1, 0, respectively). We assume that
channel state transitions are caused by primary user activity,
therefore the problem reduces to the OSA problem. For any channel
$i$, $r^i_1=1$, $r^i_0=0.1$. We simulate RCA-M in four
environments with different state transition probabilities. We
compute the normalized regret values, i.e., the regret per single
play $R(n)/M$ by averaging the results of 100 runs. 

\rev{The state} transition probabilities are given in Table \ref{tab:table1} and
the mean rewards of the channels under these state transition
probabilities are given in Table \ref{tab:table2}. 
\rev{The four environment, denoted as S1, S2, S3 and S4, respectively, are summarized
as follows. } 
In S1 channels
are bursty with mean rewards not close to each other; in S2
channels are non-bursty with mean rewards not close to each other;
in S3 there are bursty and non-bursty channels with mean rewards
not close to each other; and in S4 there are bursty and non-bursty
channels with mean rewards close to each other. 

In Figures
\ref{fig:S1_L7200}, \ref{fig:S2_L360}, \ref{fig:S3_L3600},
\ref{fig:S4_L7200}, we observe the normalized regret of RCA-M for
the minimum values of $L$ such that the logarithmic bound hold.
However, comparing with Figures \ref{fig:S1_L1}, \ref{fig:S2_L1},
\ref{fig:S3_L1}, \ref{fig:S4_L1} we see that the normalized regret
is smaller for $L=1$. Therefore the condition on $L$ we have for
the logarithmic bound, \rev{while sufficient, does not appear necessary}.  
We also observe that for the Gilbert-Elliot channel model the regret can be smaller
when $L$ is set to a value smaller than $112 S^2_{\max} r^2_{\max}
\hat{\pi}^2_{\max} /\epsilon_{\min}$.


\begin{table}[ht] 
\begin{center}
    \begin{tabular}{ l|c|c|c|c|c|c|c|c|c|c}
    \hline
    channel & 1 & 2 & 3 & 4 & 5 & 6 & 7 & 8 & 9 & 10 \\ \hline
        S1, $p_{01}$ & 0.01 & 0.01 & 0.02 & 0.02 & 0.03 & 0.03 & 0.04 & 0.04 & 0.05 & 0.05 \\ \hline
        S1, $p_{10}$ & 0.08 & 0.07 & 0.08 & 0.07 & 0.08 & 0.07 & 0.02 & 0.01 & 0.02 & 0.01 \\ \hline
        S2, $p_{01}$ & 0.1 & 0.1 & 0.2 & 0.3 & 0.4 & 0.5 & 0.6 & 0.7 & 0.8 & 0.9 \\ \hline
        S2, $p_{10}$ & 0.9 & 0.9 & 0.8 & 0.7 & 0.6 & 0.5 & 0.4 & 0.3 & 0.2 & 0.1 \\ \hline
        S3, $p_{01}$ & 0.01 & 0.1 & 0.02 & 0.3 & 0.04 & 0.5 & 0.06 & 0.7 & 0.08 & 0.9 \\ \hline
        S3, $p_{10}$ & 0.09 & 0.9 & 0.08 & 0.7 & 0.06 & 0.5 & 0.04 & 0.3 & 0.02 & 0.1 \\ \hline
        S4, $p_{01}$ & 0.02 & 0.04 & 0.04 & 0.5 & 0.06 & 0.05 & 0.7 & 0.8 & 0.9 & 0.9 \\ \hline
        S4, $p_{10}$ & 0.03 & 0.03 & 0.04 & 0.4 & 0.05 & 0.06 & 0.6 & 0.7 & 0.8 & 0.9 \\ \hline
    \end{tabular}
\caption{Transition probabilities \label{tab:table1}}
\end{center}
\end{table}

\begin{table}[ht] 
\begin{center}
    \begin{tabular}{ l|c|c|c|c|c|c|c|c|c|c}
    \hline
    channel & 1 & 2 & 3 & 4 & 5 & 6 & 7 & 8 & 9 & 10 \\ \hline
        S1 & 0.20 & 0.21 & 0.28 & 0.30 & 0.35 & 0.37 & 0.70 & 0.82 & 0.74 & 0.85 \\ \hline
        S2 & 0.19 & 0.19 & 0.28 & 0.37 & 0.46 & 0.55 & 0.64 & 0.73 & 0.82 & 0.91 \\ \hline
        S3 & 0.19 & 0.19 & 0.28 & 0.37 & 0.46 & 0.55 & 0.64 & 0.73 & 0.82 & 0.91 \\ \hline
        S4 & 0.460 & 0.614 & 0.550 & 0.600 & 0.591 & 0.509 & 0.585 & 0.580 & 0.577 & 0.550 \\ \hline
    \end{tabular}
\caption{Mean rewards \label{tab:table2}}
\end{center}
\end{table}

\begin{figure}[ht]
\begin{minipage}[b]{0.5\linewidth}
\centering
\includegraphics[scale=0.4]{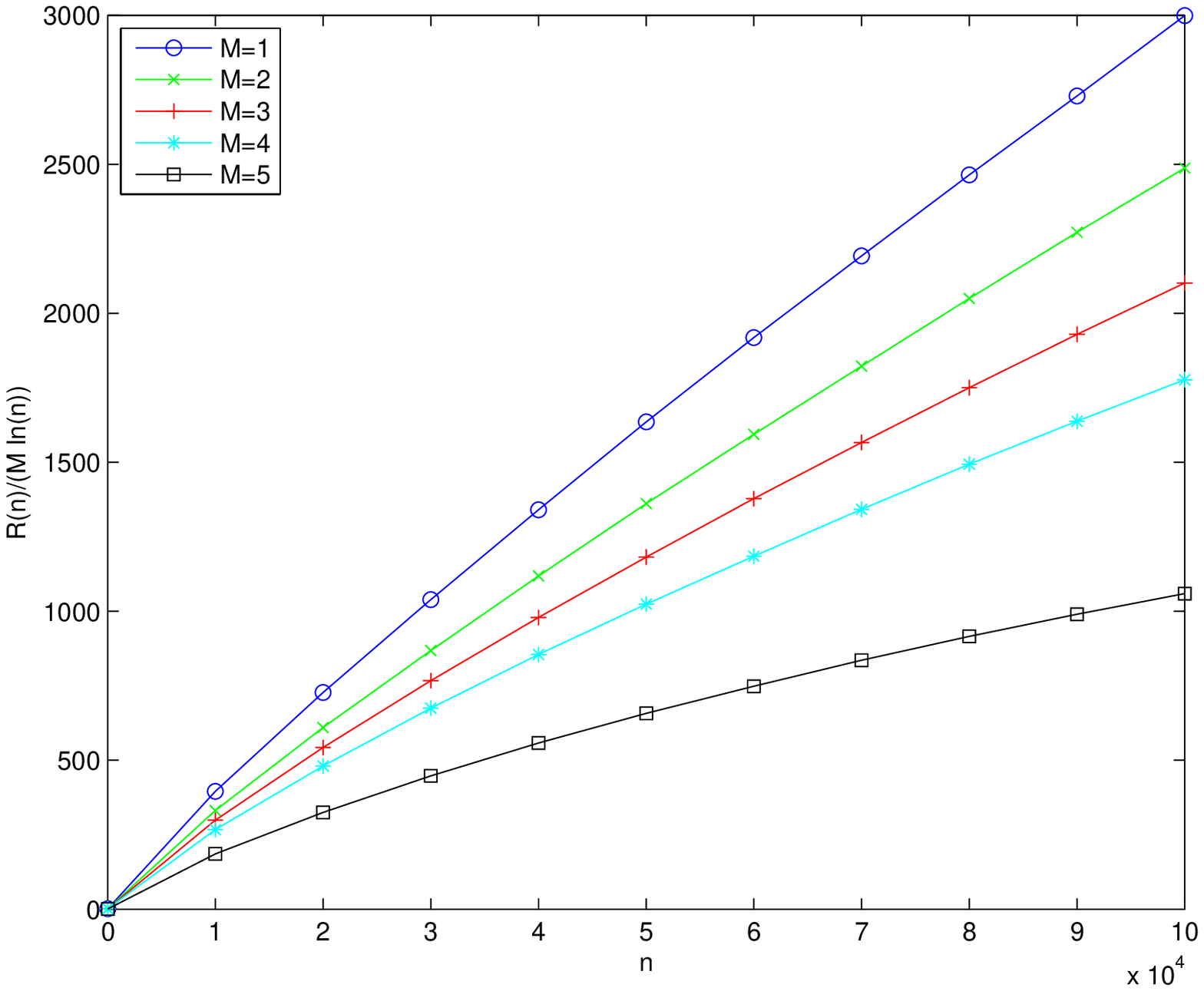}
\caption{Normalized regret under S1, $L=7200$}
\label{fig:S1_L7200}
\end{minipage}
\hspace{0.5cm}
\begin{minipage}[b]{0.5\linewidth}
\centering
\includegraphics[scale=0.4]{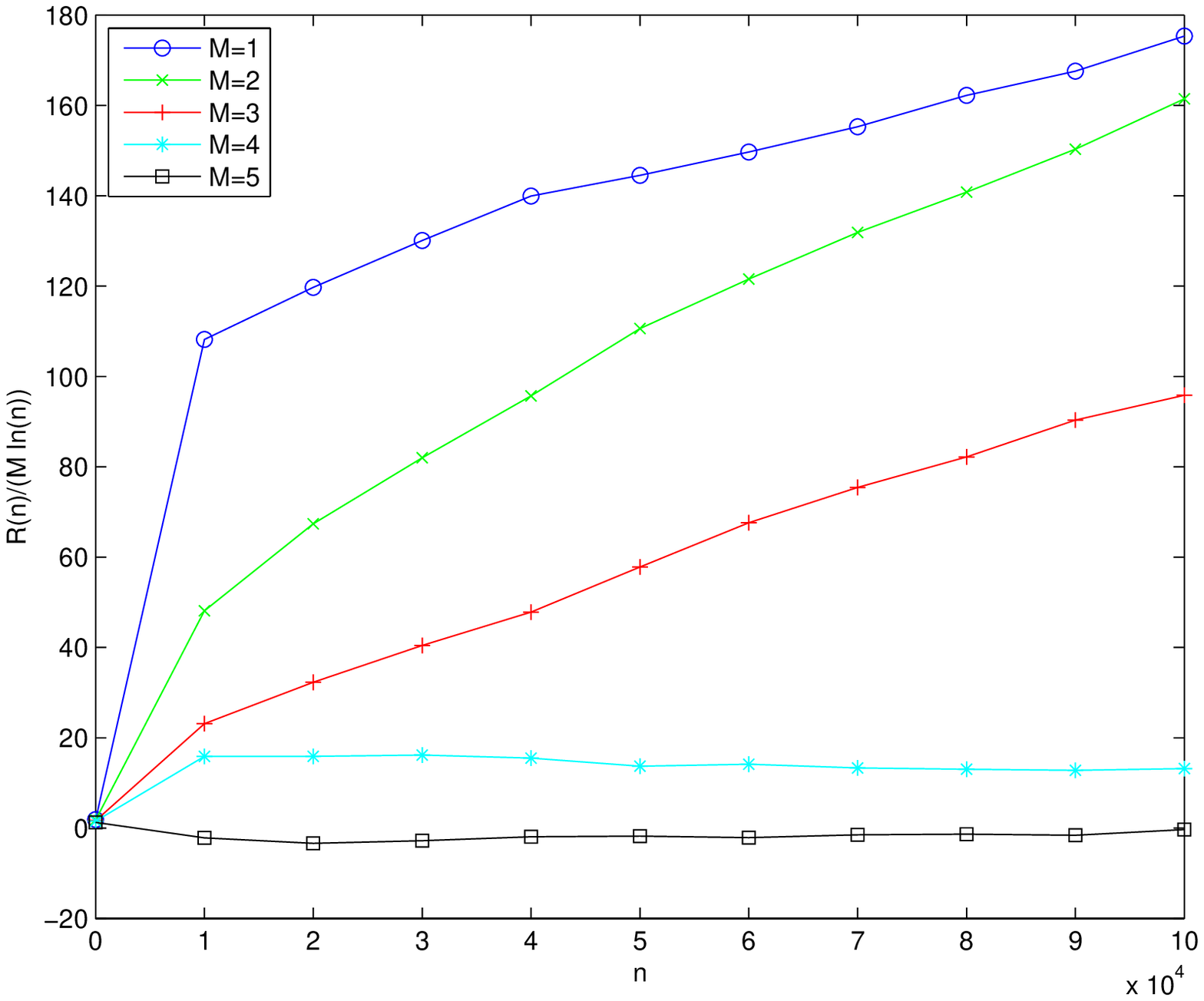}
\caption{Normalized regret under S1, $L=1$}
\label{fig:S1_L1}
\end{minipage}
\end{figure}

\comment{
\begin{figure}[ht]
\begin{center}
\includegraphics[width=3in]{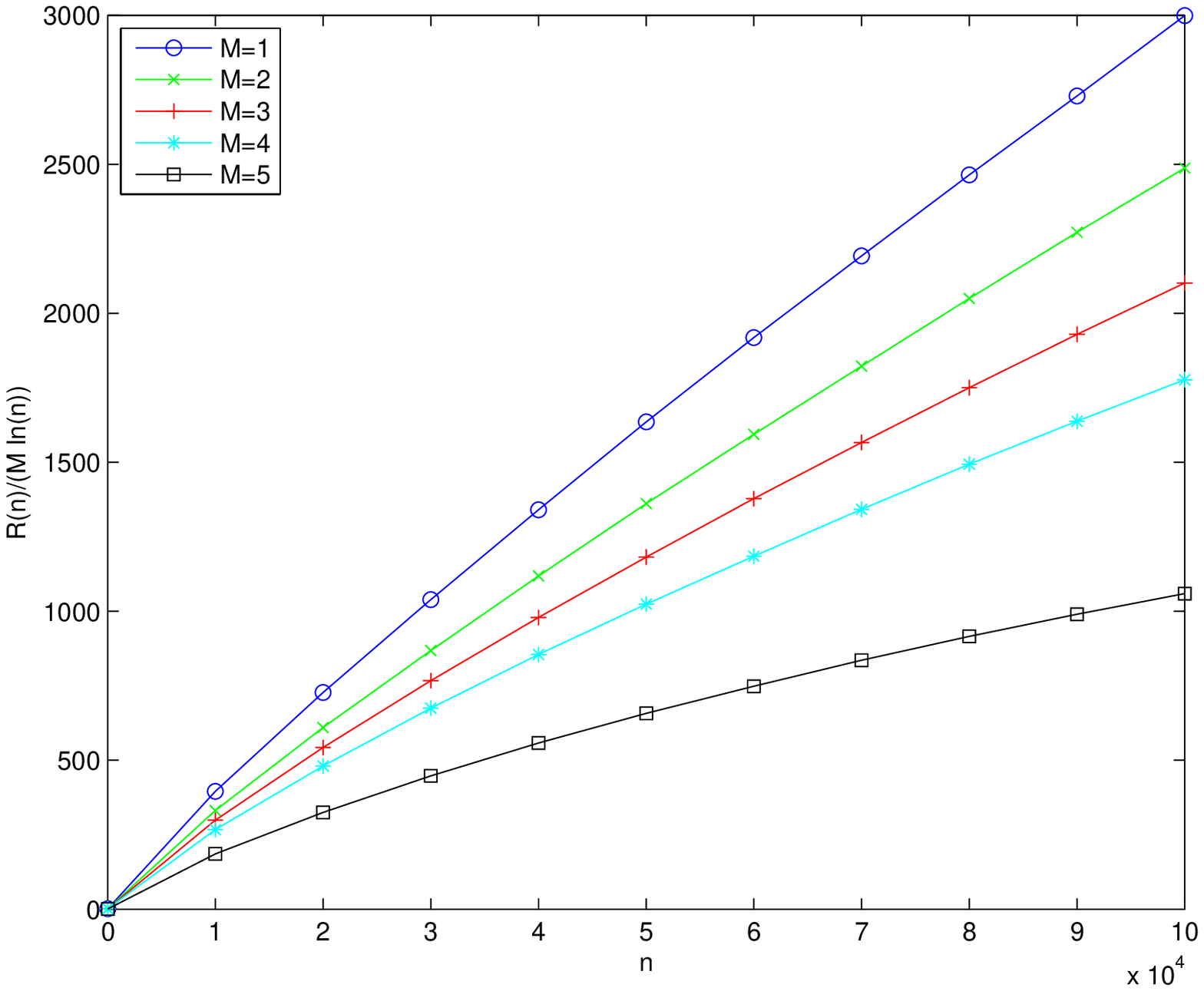}
\caption{Normalized regret under S1, $L=7200$}
\label{fig:S1_L7200}
\end{center}
\end{figure}

\begin{figure}[ht]
\begin{center}
\includegraphics[width=3in]{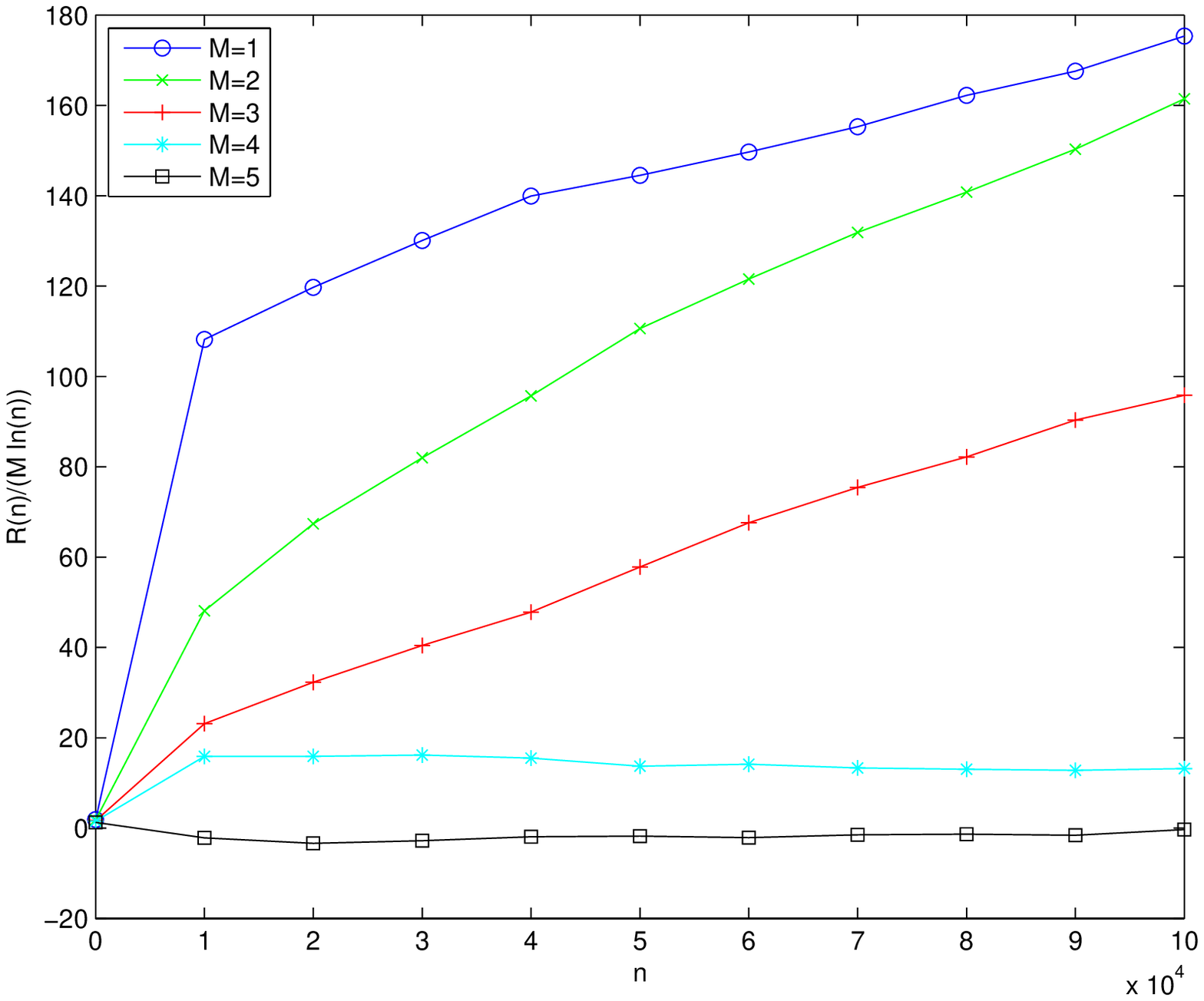}
\caption{Normalized regret under S1, $L=1$}
\label{fig:S1_L1}
\end{center}
\end{figure}
}

\begin{figure}[ht]
\begin{minipage}[b]{0.5\linewidth}
\centering
\includegraphics[scale=0.4]{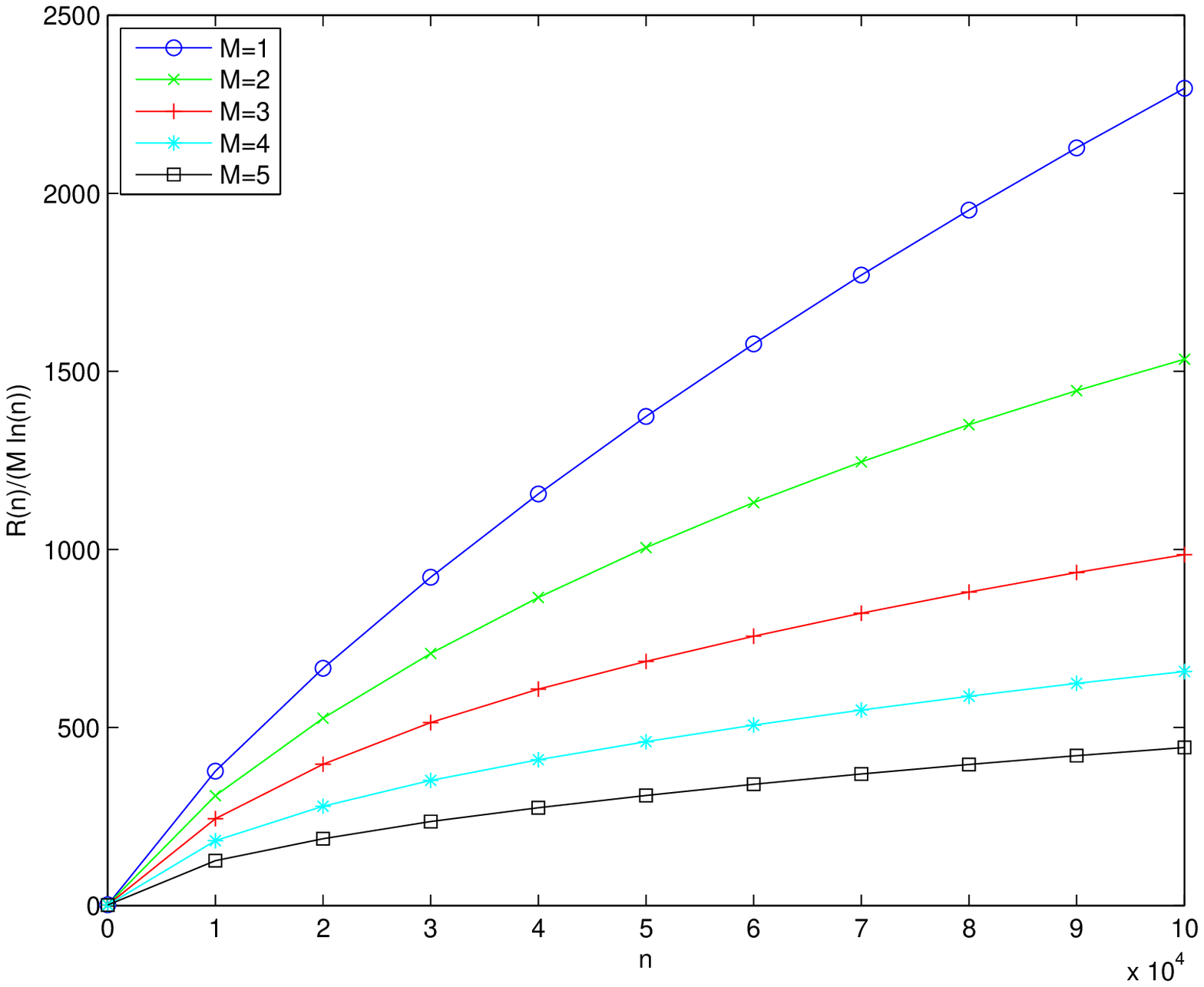}
\caption{Normalized regret under S2, $L=360$}
\label{fig:S2_L360}
\end{minipage}
\hspace{0.5cm}
\begin{minipage}[b]{0.5\linewidth}
\centering
\includegraphics[scale=0.4]{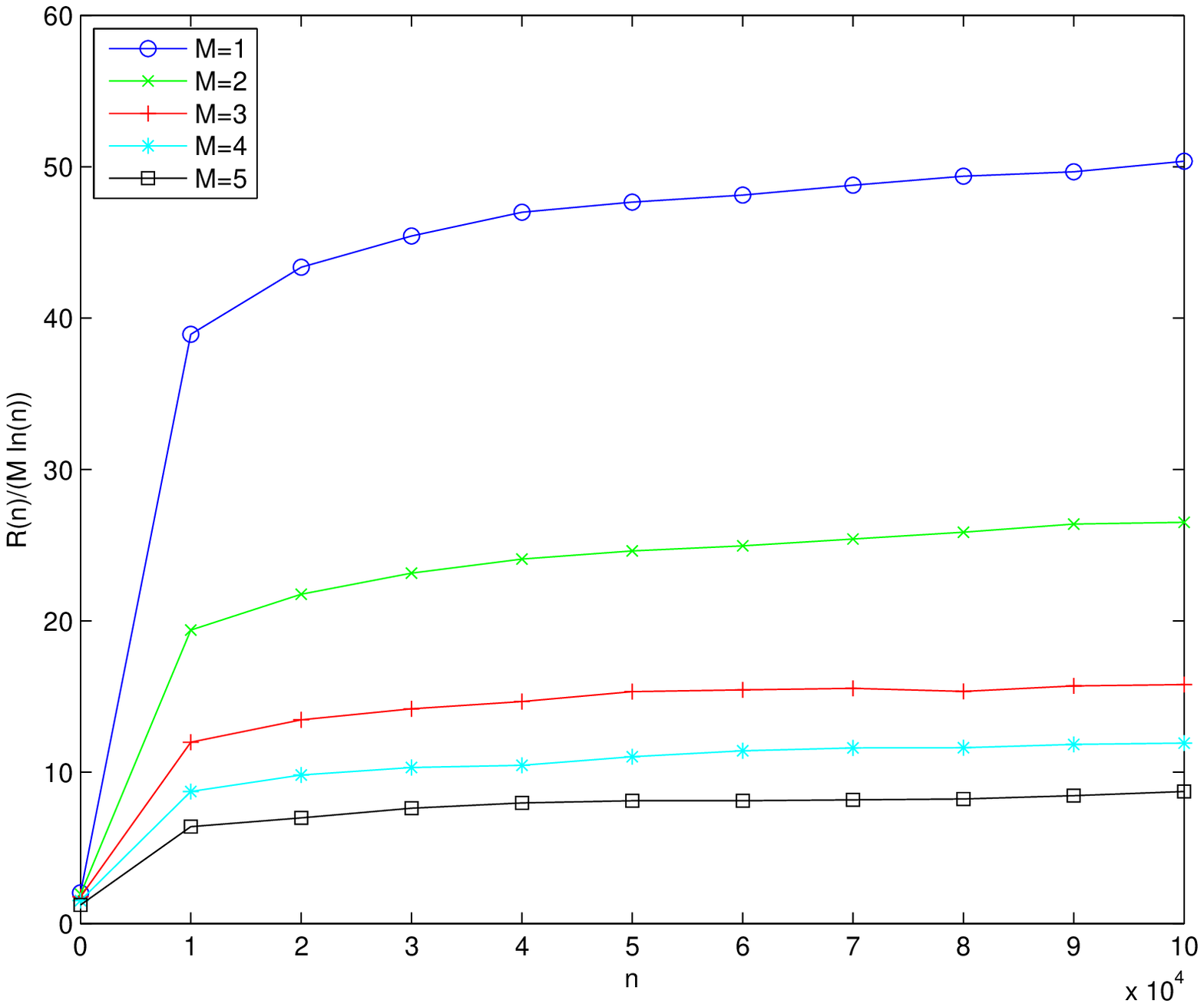}
\caption{Normalized regret under S2, $L=1$}
\label{fig:S2_L1}
\end{minipage}
\end{figure}

\comment{
\begin{figure}[ht]
\begin{center}
\includegraphics[width=3in]{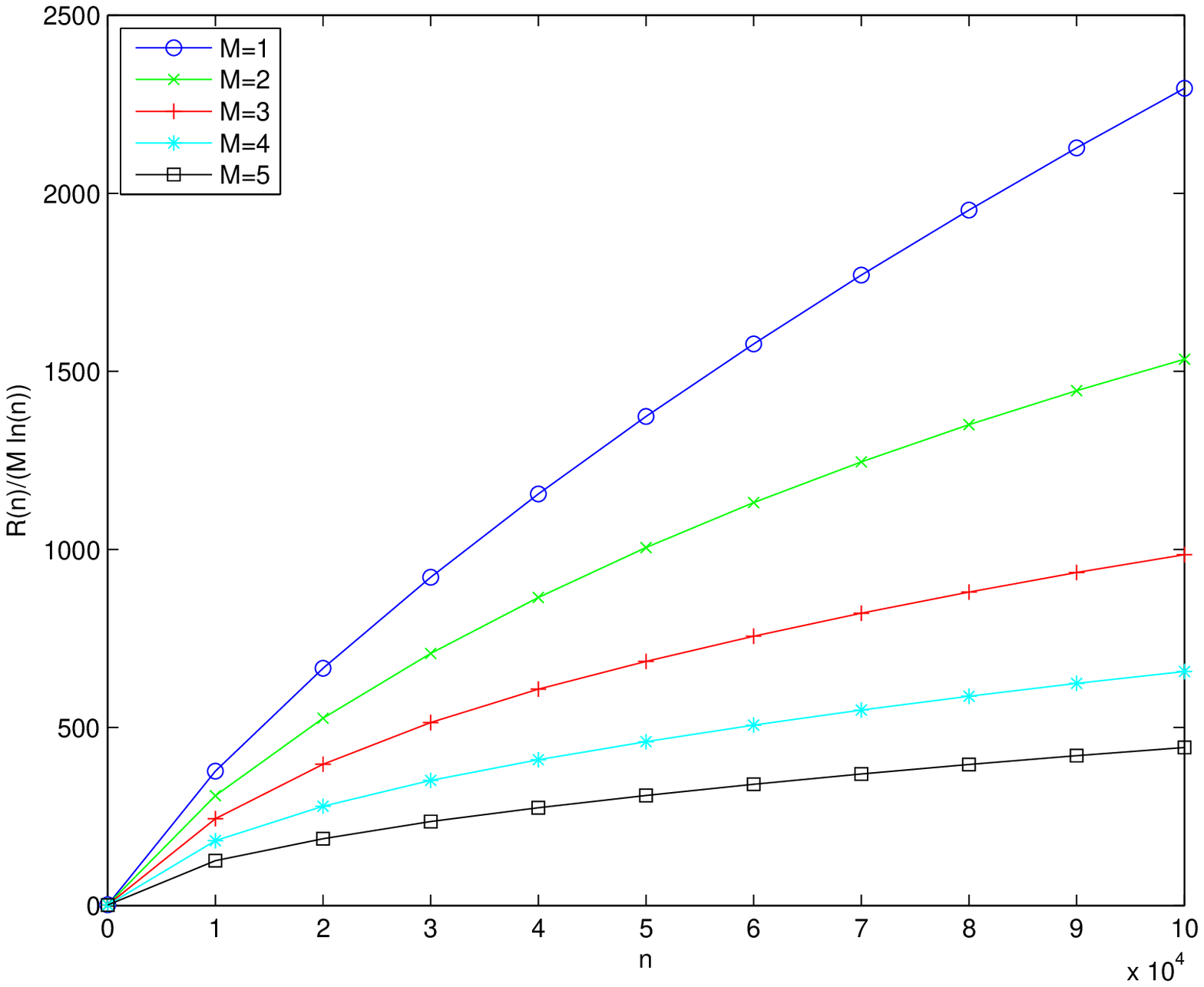}
\caption{Normalized regret under S2, $L=360$}
\label{fig:S2_L360}
\end{center}
\end{figure}

\begin{figure}[ht]
\begin{center}
\includegraphics[width=3in]{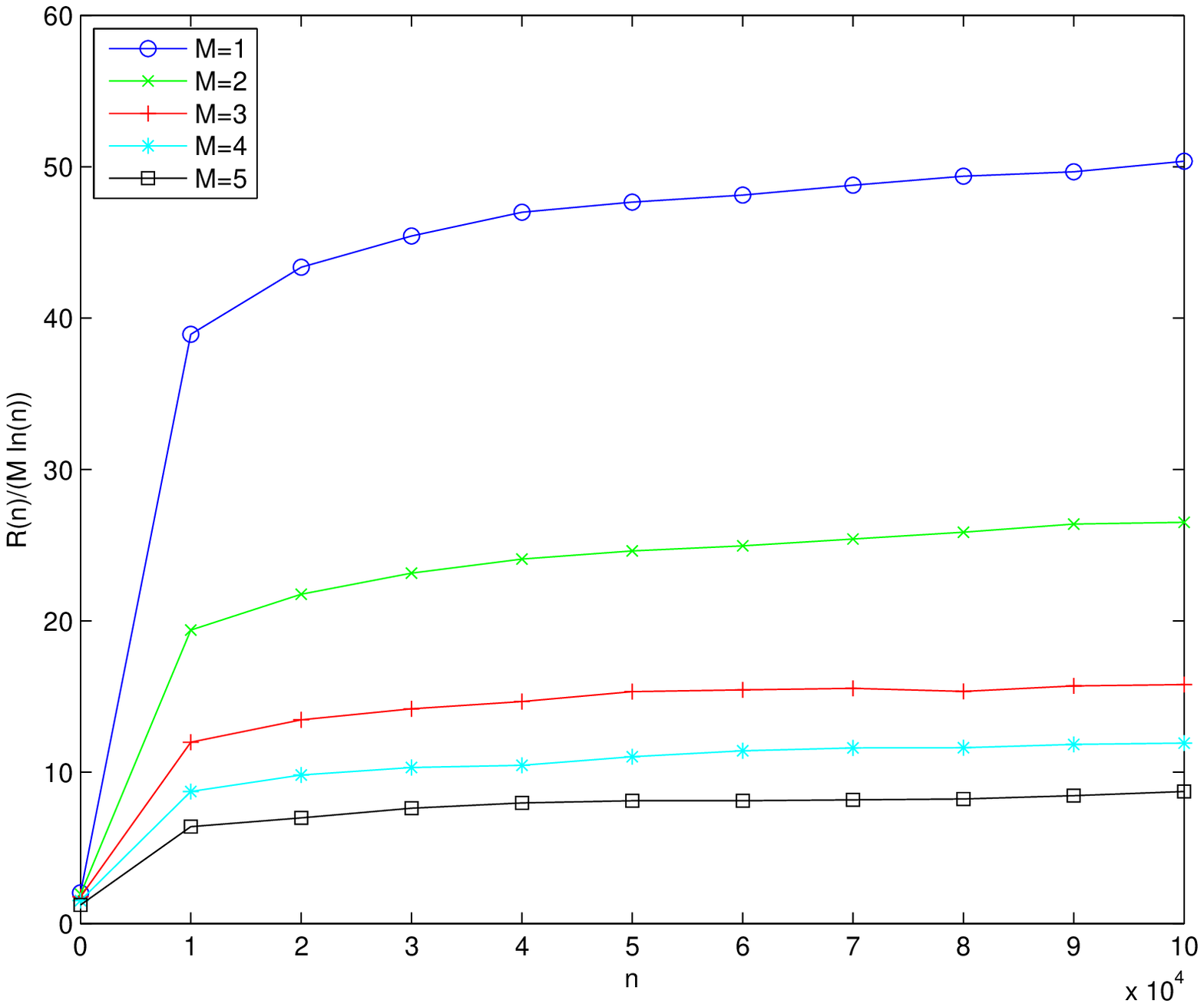}
\caption{Normalized regret under S2, $L=1$}
\label{fig:S2_L1}
\end{center}
\end{figure}
}

\begin{figure}[hb]
\begin{minipage}[b]{0.5\linewidth}
\centering
\includegraphics[scale=0.4]{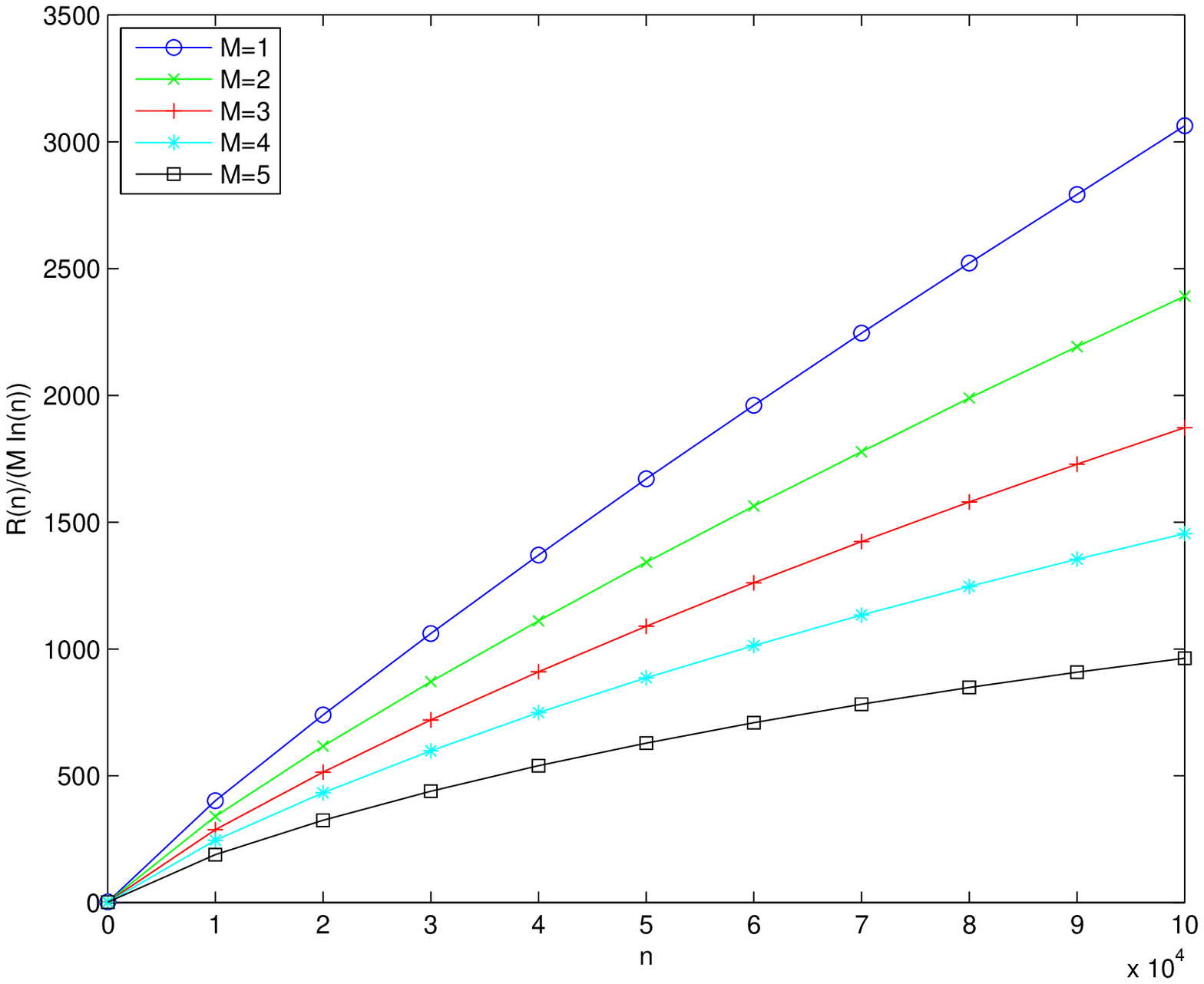}
\caption{Normalized regret under S3, $L=3600$}
\label{fig:S3_L3600}
\end{minipage}
\hspace{0.5cm}
\begin{minipage}[b]{0.5\linewidth}
\centering
\includegraphics[scale=0.4]{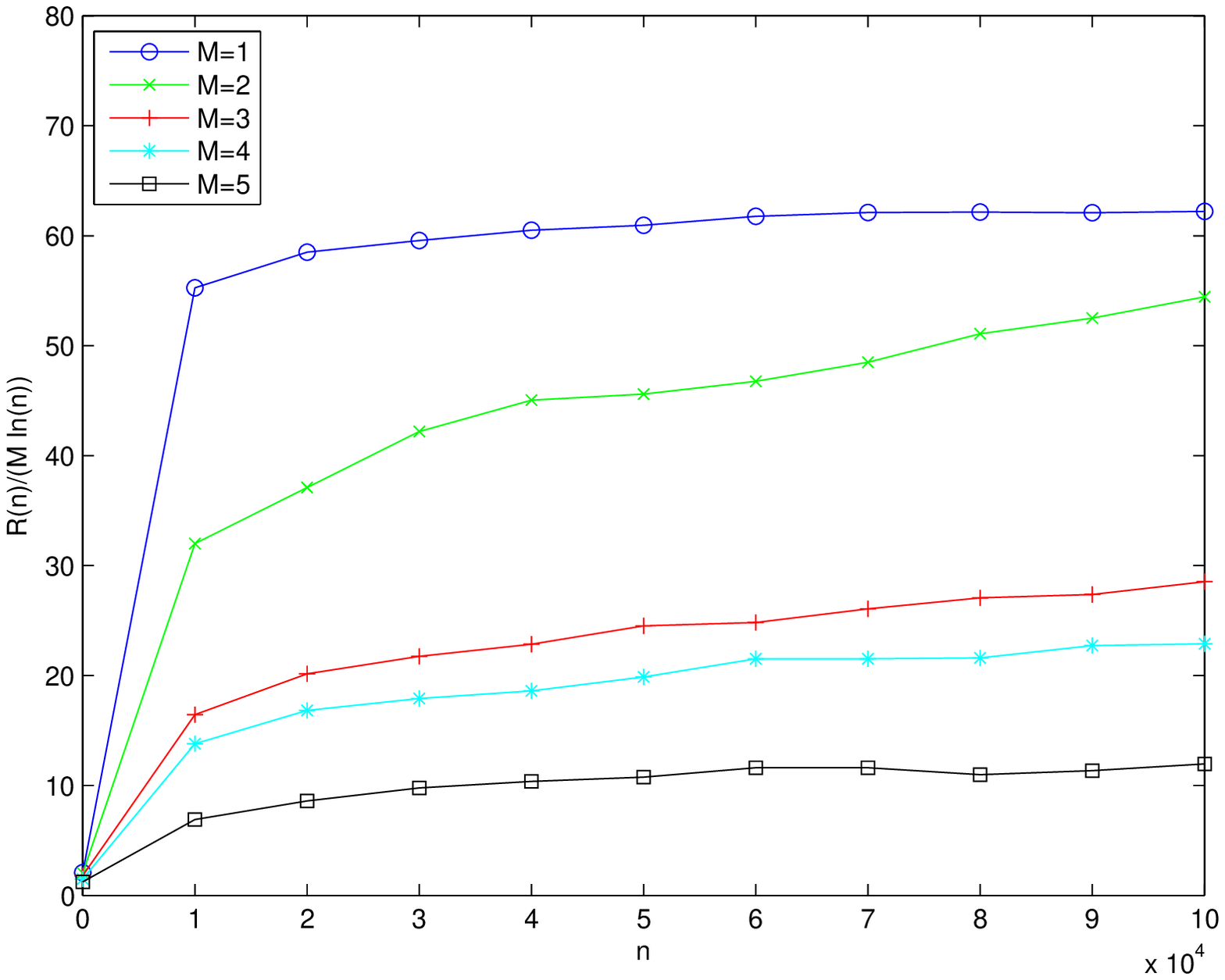}
\caption{Normalized regret under S3, $L=1$}
\label{fig:S3_L1}
\end{minipage}
\end{figure}

\comment{
\begin{figure}[ht]
\begin{center}
\includegraphics[width=3in]{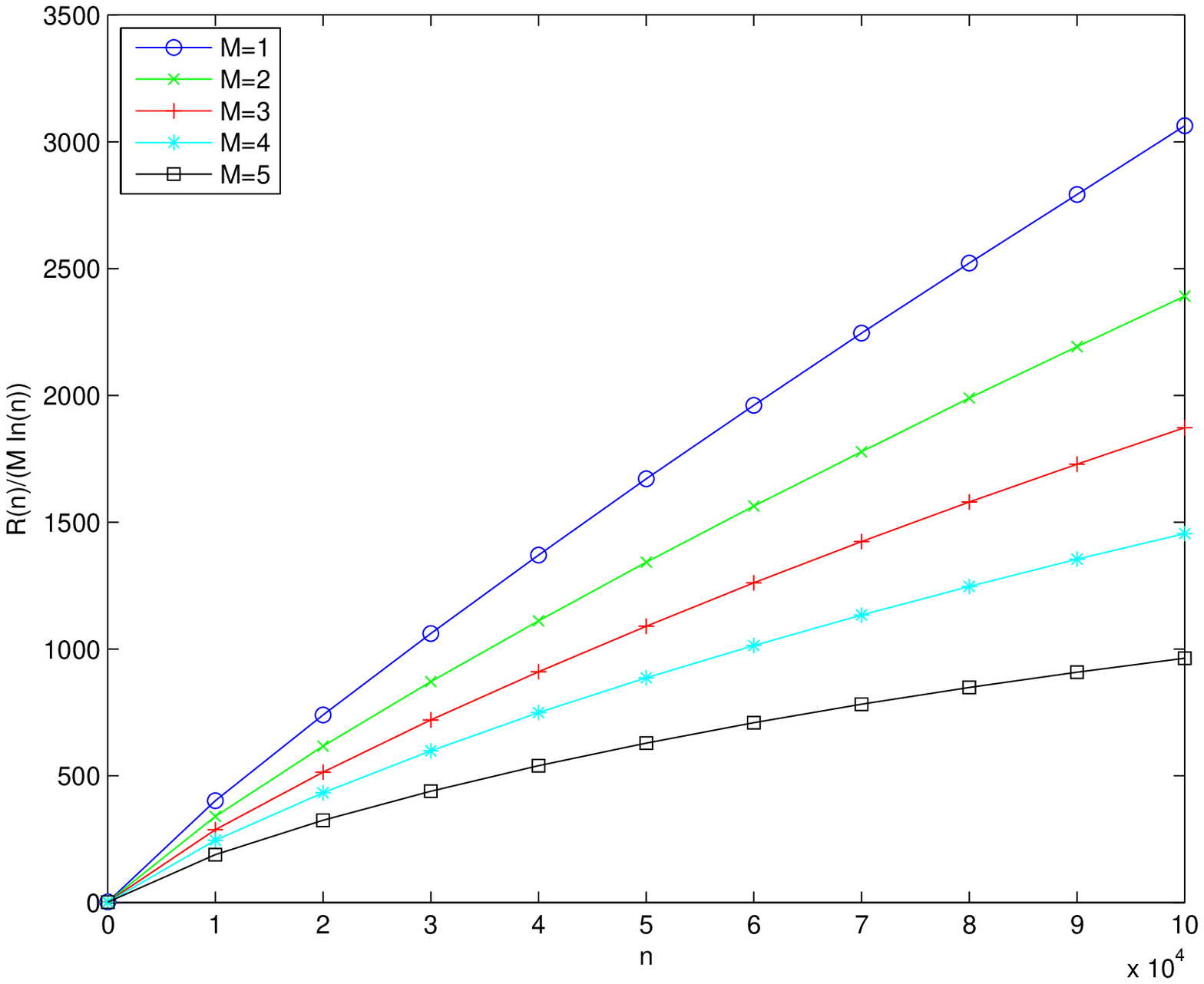}
\caption{Normalized regret under S3, $L=3600$}
\label{fig:S3_L3600}
\end{center}
\end{figure}

\begin{figure}[ht]
\begin{center}
\includegraphics[width=3in]{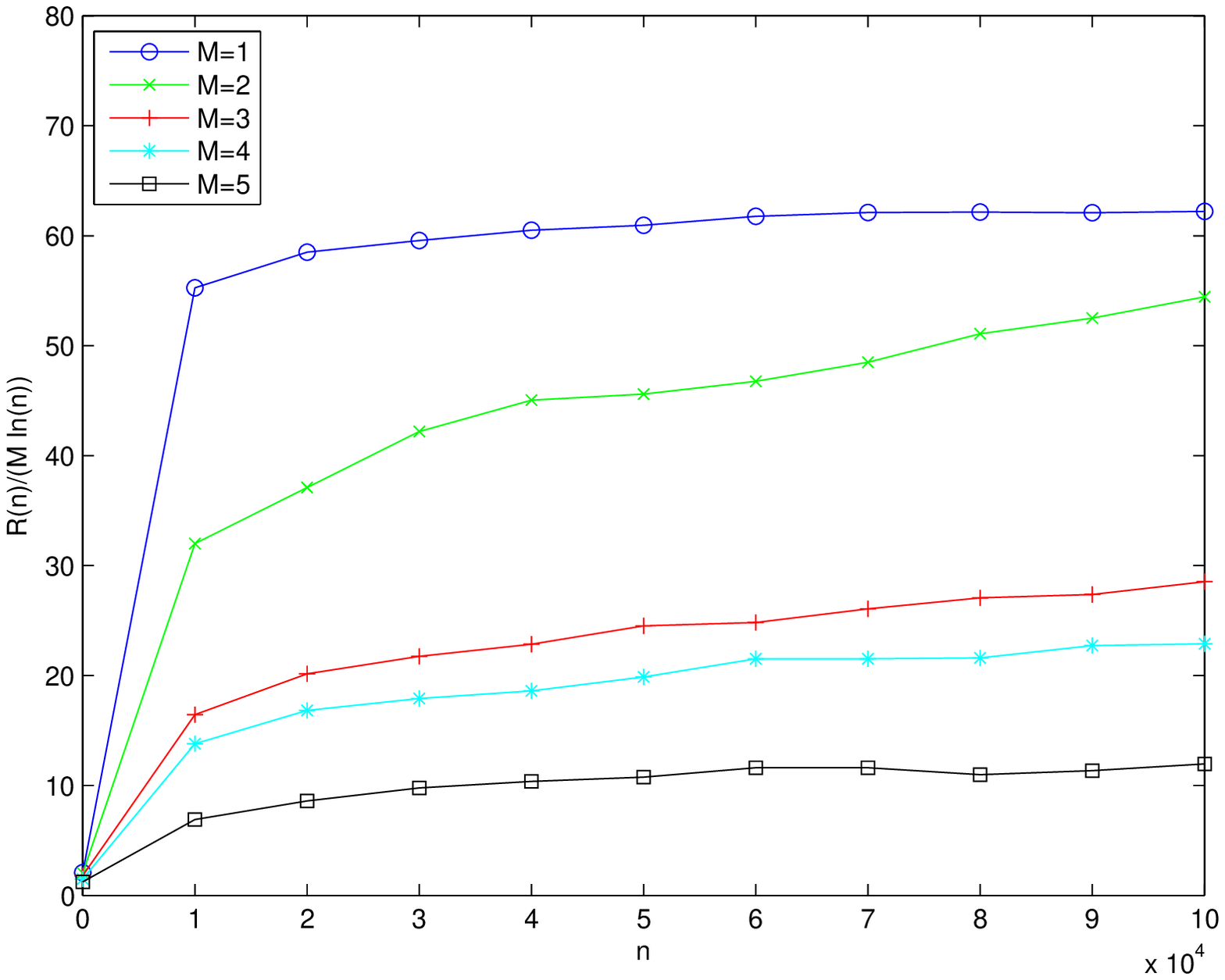}
\caption{Normalized regret under S3, $L=1$}
\label{fig:S3_L1}
\end{center}
\end{figure}
}

\begin{figure}[h]
\begin{minipage}[b]{0.5\linewidth}
\centering
\includegraphics[scale=0.4]{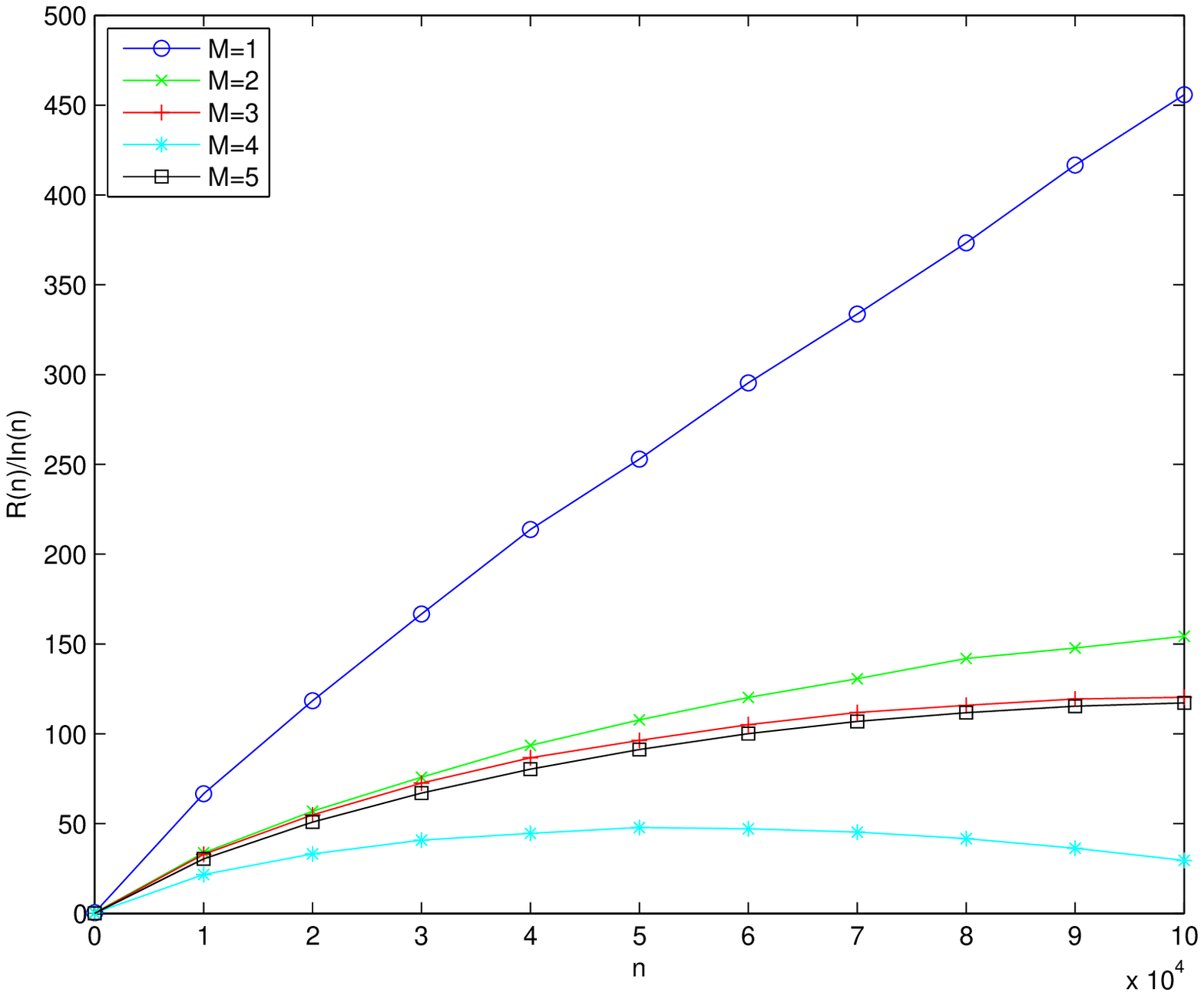}
\caption{Normalized regret under S4, $L=7200$}
\label{fig:S4_L7200}
\end{minipage}
\hspace{0.5cm}
\begin{minipage}[b]{0.5\linewidth}
\centering
\includegraphics[scale=0.4]{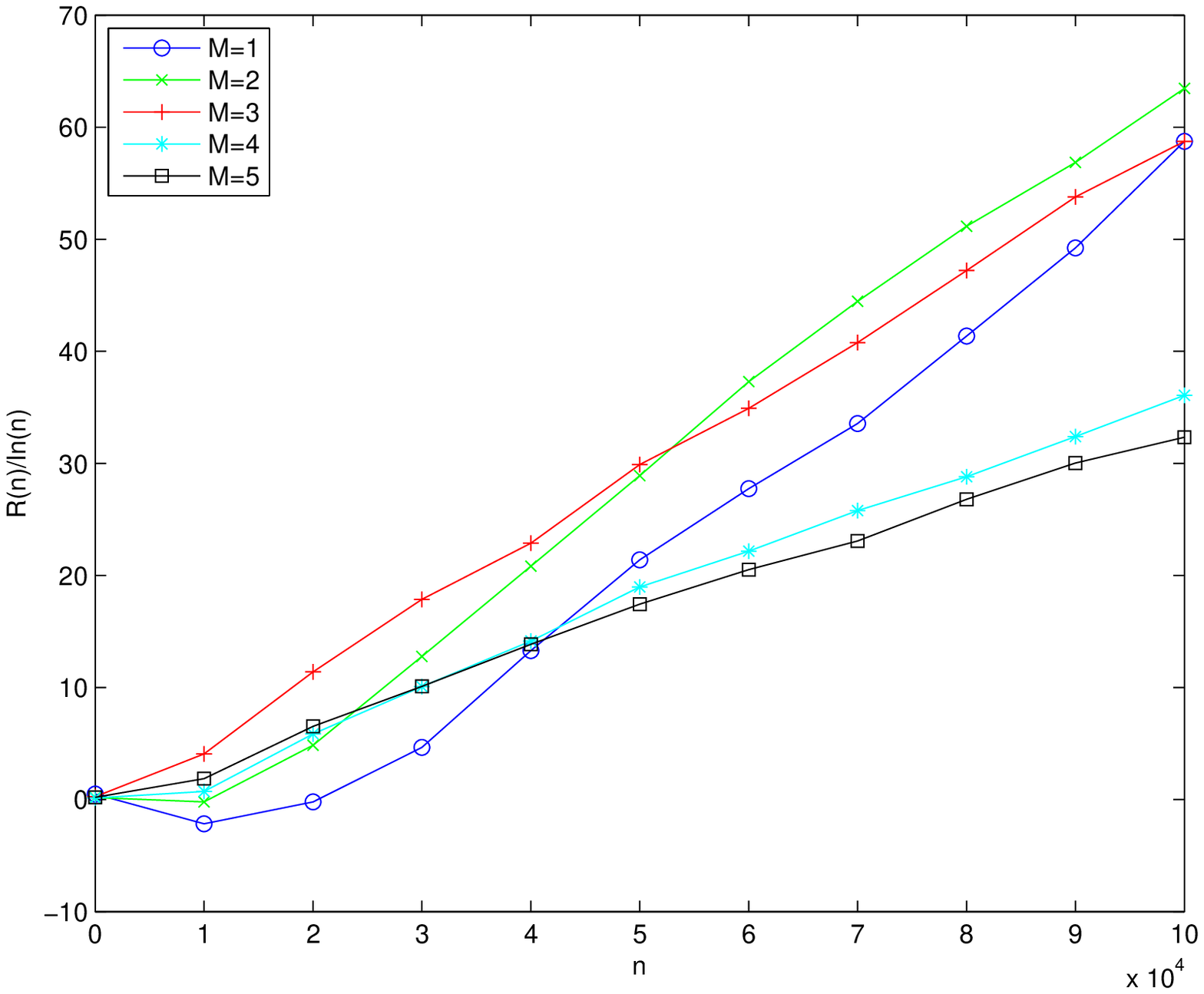}
\caption{Normalized regret under S4, $L=1$}
\label{fig:S4_L1}
\end{minipage}
\end{figure}

\comment{
\begin{figure}[ht]
\begin{center}
\includegraphics[width=3in]{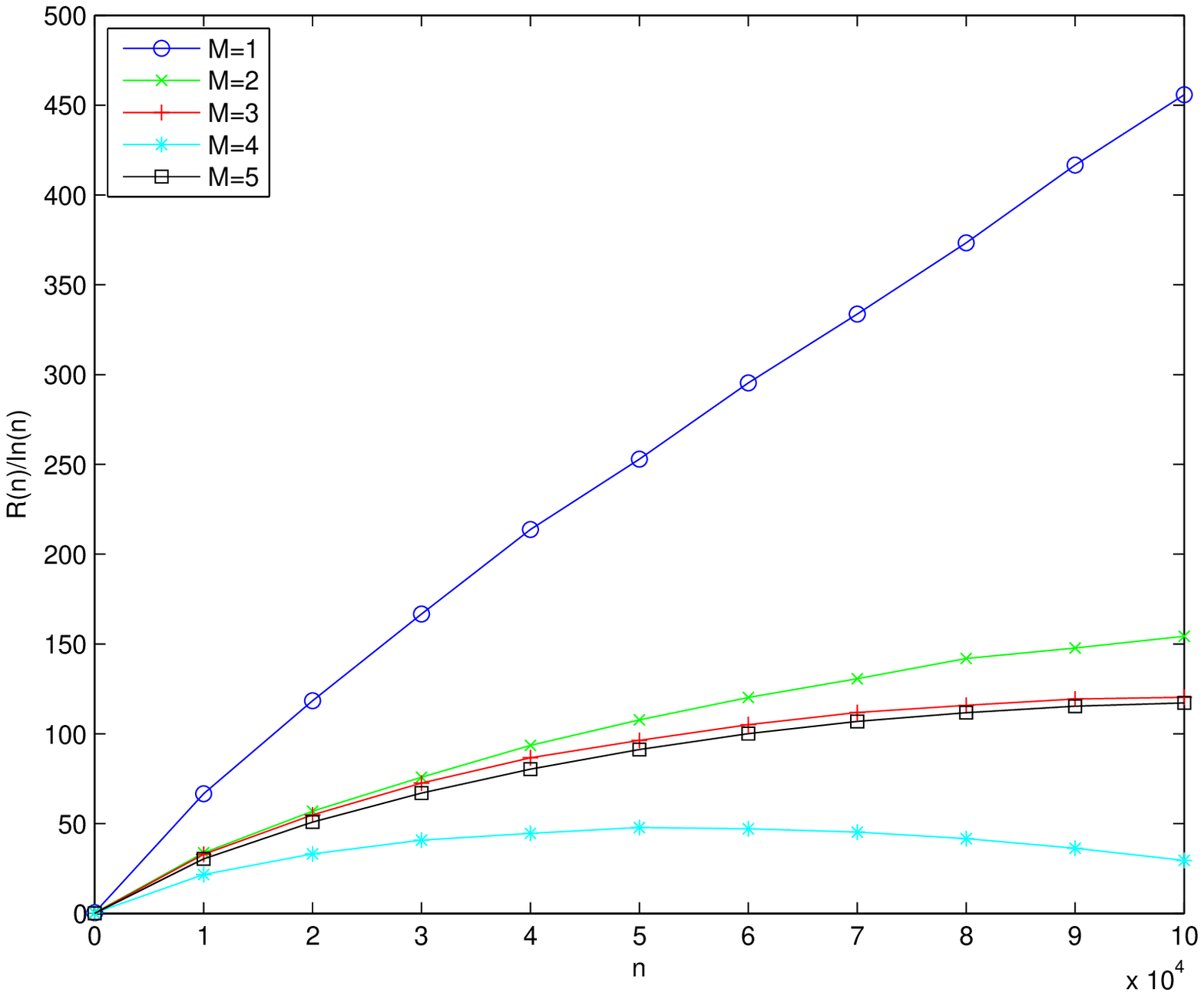}
\caption{Normalized regret under S4, $L=7200$}
\label{fig:S4_L7200}
\end{center}
\end{figure}

\begin{figure}[ht]
\begin{center}
\includegraphics[width=3in]{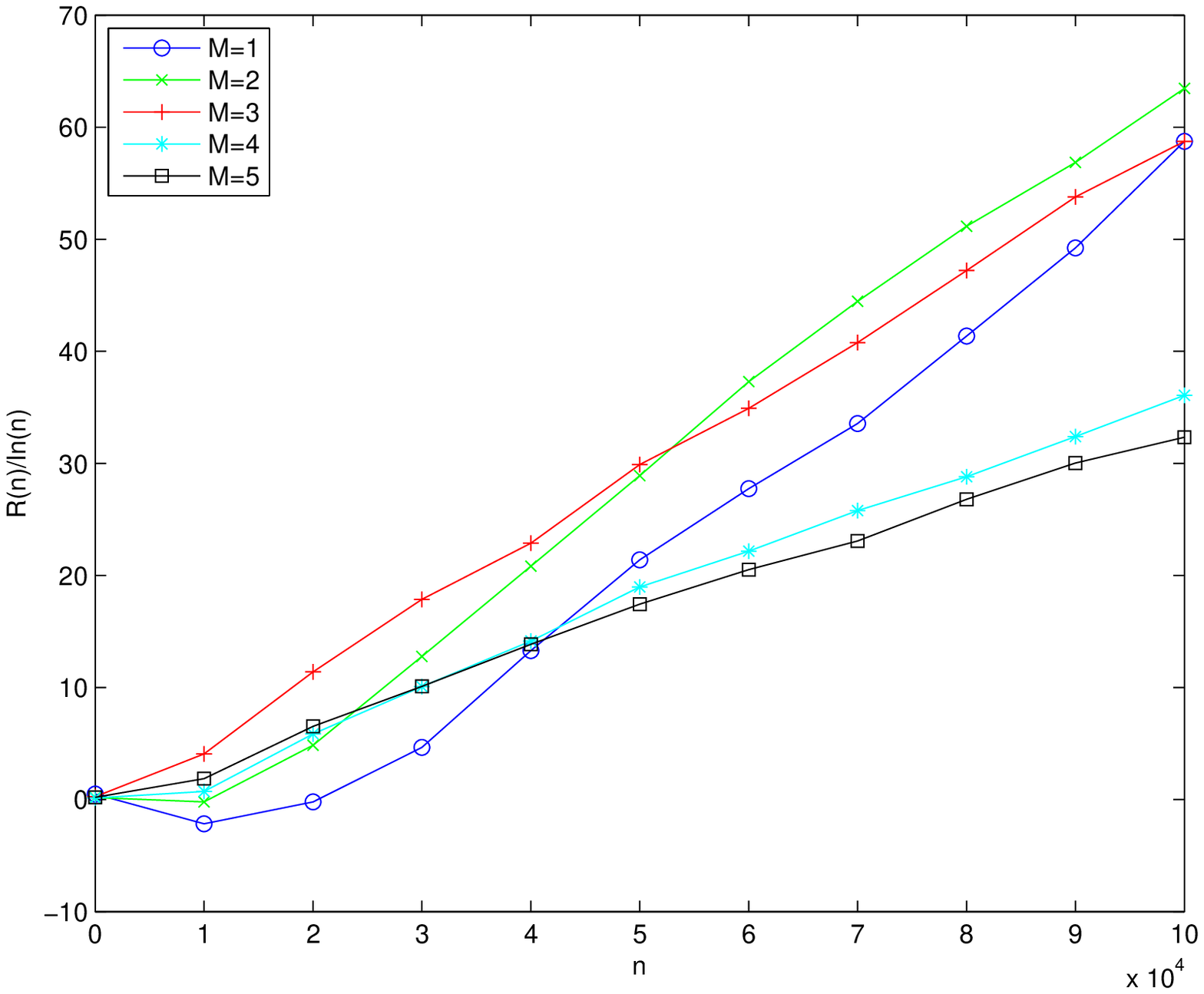}
\caption{Normalized regret under S4, $L$=1}
\label{fig:S4_L1}
\end{center}
\end{figure}
}
\section{Discussion} \label{sec:discussion} 

\rev{In this section we discuss how the performance of RCA-M may be improved (in terms of the constants and not in order), and possible relaxation and extensions.} 

\subsection{Applicability, Performance Improvement, and Relaxation}

We note that the same logarithmic bound \rev{derived in this paper holds for the general restless bandit where the state evolution is given by two matrices: the active and passive transition probability matrices ($P^i$ and $Q^i$ respectively for arm $i$), which are potentially different.  The addition of a different $Q^i$ does not affect the analysis because} the reward to the player from an arm is determined only by the active transition probability matrix and the first state after a discontinuity in playing the arm. Since the number of plays from any suboptimal arm is logarithmic and the expected hitting time of any state is finite the regret due to $Q^i$ is at most logarithmic. 
\rev{We further note that for the same reason the arm may not even follow a Markovian rule in the passive state, and the same logarithmic bound will continue to hold.} 

The regenerative state \rev{for an arm under} RCA-M \rev{is} chosen based on the
\rev{random} initial observation.  \rev{This means that RCA-M may happen upon} a state with long
recurrence time which will result in long SB1 and SB2 \rev{sub-blocks}. We \rev{propose the following modification}: RCA-M \rev{records} all observations from all arms. Let $k_i(s,t)$ be the \rev{total} number
of observations from arm $i$ \rev{up to} time $t$ \rev{that are {\em excluded}} from \rev{the} computation of \rev{the} index \rev{of} arm $i$ \rev{when the regenerative state is $s$}.  \rev{Recall that the index of an arm is computed based on observations from regenerative cycles; this implies that $k_i(s,t)$ is the total number of slots in SB1's when the regenerative state is $s$.} 
Let $t_n$ be the time at the end of the
$n$-th block. If the arm to be played in the $n$-th block is $i$ then
the regenerative state is set to $\gamma^i(n) = \arg\min_{s \in
S^i} k_i(s,t_{n-1})$. The idea behind this modification is to
estimate the state with the smallest recurrence time and choose
the regenerative cycles according to this state. With this
modification the number of observations that does not contribute
to the index computation and the probability of choosing a
suboptimal arm \rev{can be minimized over time}. 

It's also worth noting that the selection of the regenerative state $\gamma^i$ in each block in general can be arbitrary: within the same SB2, we can start and end in different states.  As long as we guarantee that two successive SB2's end and start with the same state, we will have a continuous sample path for which our analysis in Section \ref{sec:restless-problem} holds.

\subsection{Relaxation of Certain Conditions} 

\rev{We have noted in Section \ref{sec:example} that the condition on $L$ while sufficient does not appear necessary for the logarithmic regret bound to hold.  Indeed our examples show} that smaller regret can be achieved by setting $L=1$.  
\rev{Note that this condition on $L$ originates from the large deviation bound by Lezaud given in Lemma \ref{lemma3}.}  This condition can be relaxed if we use a tighter large deviation bound. 

\rev{We further note that even if no information is available on the underlying Markov chains to derive this sufficient condition on $L$,
an $o(log(n)f(n))$ regret is achievable by letting $L$ grow slowly with time where $f(n)$ is any increasing sequence. 
Such approach has been used in other settings and algorithms, see e.g., \cite{anandkumar, hliu1}.} 

\rev{We have noted earlier that the strict inequality $\mu^M > \mu^{M+1}$ is required for the restless multiarmed bandit problem because in order to have logarithmic regret, we can have no more than a logarithmic number of discontinuities from the optimal arms}. When $\mu^M = \mu^{M+1}$ the rankings of the indices of arms $M$ and $M+1$ can oscillate indefinitely resulting in a large number of discontinuities. 
\rev{Below we briefly discuss how to resolve this issue if indeed $\mu^M = \mu^{M+1}$.}  Consider adding a threshold $\epsilon$ to the algorithm such that a new arm will be selected instead of an arm currently being played only if the index of that arm is at least $\epsilon$ larger than the index of the currently played arm which has the smallest index among all currently played arms. Then given that $\epsilon$ is sufficiently small (with respect to the differences of mean rewards) 
indefinite switching between the $M$-th and the $M+1$-th arms can be avoided.  
\rev{However, further analysis is needed to verify that this approach will result in logarithmic regret.} 

\subsection{Definition of Regret}

We have used the weak regret measure throughout this paper, which compares the learning strategy with the best single-action strategy.  When the statistics are known a priori, it is clear that in general the best one can do is not a single-action policy (in principle one can drive such a policy using dynamic programming).  Ideally one could try to adopt a regret measure with respect to this optimal policy.  However, such an optimal policy in the restless case is not known in general \cite{whittle,papadim}, which makes the comparison intractable, except for some very limited cases when such a policy happens to be known \cite{ahmad, wdai1}. 





\subsection{Extensions to A Decentralized Multiplayer Setting and Comparison with Similar Work}
\rev{As mentioned in the introduction, there has been a number of recent studies extending single player algorithms to multi-player settings where collisions are possible \cite{hliu2, anandkumar}.  Within this context we note that RCA-M in its currently form does not extend in a straightforward way to a decentralized multi-player setting.} 
It \rev{remains an interesting subject of future study.} 
A recent work \cite{hliu1} considers the \rev{same} restless multiarmed bandit
problem studied in the present paper. They achieve logarithmic regret by using
exploration and exploitation blocks that grow \rev{geometrically} with time. 
The construction in \cite{hliu1} is very different from ours, but is amenable to multi-player extension \cite{hliu2} due to the constant, though growing, nature of the block length which can be synchronized among players. 

It is interesting to note that the essence behind our approach RCA-M is to reduce a restless bandit problem to a rested bandit problem; this done by sampling in a way to construct a continuous sample path, which then allows us to use the same set of large deviation bounds over this reconstructed, entire sample path.  By contrast, the method introduced in \cite{hliu1} applies large deviation bounds to individual segments (blocks) of the observed sample path (which is not a continuous sample path representative of the underlying Markov chain because the chain is restless); this necessitates the need to precisely control the length and the number of these blocks, i.e., they must grow in length over time. 
\rev{
Another difference is that under our scheme, the exploration and exploitation are done simultaneously and implicitly through the use of the index, whereas under the scheme in \cite{hliu1}, the two are done separately and explicitly through two different types of blocks. } 

\section{Conclusion} \label{sec:conc}
\rev{
In this paper we considered the rested and restless multiarmed bandit problem with Markovian rewards and multiple plays.  We showed that a simple extension to UCB1 produces logarithmic regret uniformly over time. We then constructed an algorithm RCA-M that utilizes regenerative
cycles of a Markov chain to compute a sample mean based index policy. The sampling approach reduces a restless bandit problem to the rested version, and we showed that under mild conditions on the state transition probabilities of the Markov chains this algorithm
achieves logarithmic regret uniformly over time for the restless bandit problem, and that this regret bound is also optimal. 
We numerically examine the performance of this algorithm in the case of an OSA problem with the Gilbert-Elliot
channel model.} 


\bibliographystyle{IEEE}
\bibliography{OSA}

\appendices
\section{proof of \rev{Lemma \ref{lemma:rested1}}} \label{app:B}

Let $X^{i,j}(t)$ be the state observed from the $t$th play of arm $i$ by player $j$ and $T^{i,j}(n)$ be the total number of times player $j$ played arm $i$ up to and including time $n$. Then we have,

\begin{eqnarray}
&& \left|R(n)-\left(n\sum_{j=1}^M \mu^{j}- \sum_{i=1}^K \mu^i E[T^i(n)]\right)\right| \nonumber \\ 
&=& \left|E\left[\sum_{j=1}^M \sum_{i=1}^K \sum_{x \in S^i} r^i_x \sum_{t=1}^{T^{i,j}(n)} I(X^{i,j}(t)=x)\right] - \sum_{j=1}^M \sum_{i=1}^K \sum_{x \in S^i} r^i_x \pi^i_x E[T^{i,j}(n)]\right| \nonumber \\
&=& \left| \sum_{j=1}^M \sum_{i=1}^K \sum_{x \in S^i} r^i_x (E[N^j(x, T^{i,j}(n))]- \pi^i_x E[T^{i,j}(n)])\right| \nonumber \\
&\leq& \sum_{j=1}^M \sum_{i=1}^K \sum_{x \in S^i} r^i_x C_{P^i} = C_{\mathbf{S,P,r}} \label{eqn:rested1-1} 
\end{eqnarray}
where
\begin{eqnarray}
N^j(x, T^{i,j}(n))=\sum_{t=1}^{T^{i,j}(n)} I(X^{i,j}(t)=x), \nonumber
\end{eqnarray}
and (\ref{eqn:rested1-1}) follows from Lemma \ref{lemma:anantharam} using the fact that $T^{i,j}(n)$ is a stopping time with respect to the $\sigma$-field generated by the arms played up to time $n$.

\section{} \label{app:A}
\begin{lemma}
\label{lemma-key}
%
\rev{Assume Condition \ref{cond:1} holds and all arms are rested. } 
Let $g^i_{t,s}=\bar{r}^i(s)+c_{t,s}$, $c_{t,s}=\sqrt{L\ln t/s}$.
\rev{Under} UCB-M with constant $L \geq 112 S^2_{\max}
r^2_{\max} \hat{\pi}^2_{\max} /\epsilon_{\min}$, 
for any suboptimal arm $i$ and optimal arm $j$ we have
\begin{eqnarray}
E \left[\sum_{t=1}^{n} \sum_{w=1}^{t-1} \sum_{w_i=l}^{t-1} I(g^{j}_{t, w} \leq g^i_{t, w_i})\right] \leq \frac{|S^i|+|S^j|}{\pi_{\min}} \beta,
\end{eqnarray}
where
\rev{$l=\left\lceil \frac{4L\ln n}{(\mu^M-\mu^i)^2}\right\rceil$ and} $\beta = \sum_{t=1}^{\infty} t^{-2}$.
\end{lemma}

\begin{proof}
First, we show that for any suboptimal arm $i$ and optimal arm $j$, \rev{we have that} $g^j_{t, w} \leq g^i_{t, w_i}$ implies at least one of the following holds: 
\begin{eqnarray}
\bar{r}^j(w) &\leq& \mu^j - c_{t,w} \label{eqn2}\\ 
\bar{r}^i(w_i) &\geq& \mu^i+c_{t,w_i} \label{eqn3}\\
\mu^j &<& \mu^i+2c_{t,w_i} \label{eqn4}.
\end{eqnarray}
This is because if none of the above holds, then we must have 
\begin{eqnarray*}
g^j_{t, w} = \bar{r}^j(w)+c_{t,w} > \mu^j \geq \mu^i+2c_{t,w_i} > \bar{r}^i(w_i)+c_{t,w_i} = g^i_{t, w_i}, 
\end{eqnarray*}
which contradicts $g^j_{t, w} \leq g^i_{t, w_i}$. 

If we choose $w_i \geq 4L\ln n/(\mu^{M}-\mu^i)^2$, then
\begin{eqnarray*}
2c_{t,w_i} = 2\sqrt{\frac{L \ln t}{w_i}} \leq 2\sqrt{ \frac{L \ln t (\mu^{M}-\mu^i)^2}{4L \ln n}} \leq \mu^j- \mu^i \text{ for } t \leq n, 
\end{eqnarray*}
which means (\ref{eqn4}) is false, and therefore at least one of (\ref{eqn2}) and (\ref{eqn3}) is true with this choice of $w_i$.  
%
Let $l=\left\lceil \frac{4L\ln n}{(\mu^{M}-\mu^i)^2}\right\rceil$. Then we have,
\begin{eqnarray*}
 E\left[\sum_{t=1}^{n} \sum_{w=1}^{t-1} \sum_{w_i=l}^{t-1} I(g^{j}_{t, w} \leq g^i_{t, w_i})\right] &\leq& \sum_{t=1}^{n}\sum_{w=1}^{t-1}\sum_{w_i=\left\lceil \frac{4L\ln n}{(\mu^{M}-\mu^i)^2}\right\rceil}^{t-1} \left(P(\bar{r}^j(w)\leq \mu^j-c_{t,w}) + P(\bar{r}^i(w_i) \geq \mu^i+c_{t,w_i})\right) \\
 &\leq& \sum_{t=1}^{\infty}\sum_{w=1}^{t-1}\sum_{w_i=\left\lceil \frac{4L\ln n}{(\mu^{M}-\mu^i)^2}\right\rceil}^{t-1} \left(P(\bar{r}^j(w)\leq \mu^j-c_{t,w}) + P(\bar{r}^i(w_i) \geq \mu^i+c_{t,w_i})\right) ~. 
\end{eqnarray*}

Consider an initial distribution ${\bf q}^i$ for the $i$th arm.  We have: 
\begin{eqnarray*}
&& N_{\mathbf{q}^i}=\left\|\left(\frac{q_{y}^i}{\pi_{y}^i},y\in S^i\right)\right\|_2 \leq \sum_{y \in S^i} \left\|\frac{q_{y}^i}{\pi_{y}^i}\right\|_2 \leq \frac{1}{\pi_{\min}} , 
\end{eqnarray*}
where the first inequality follows from the Minkowski inequality.
Let $n^i_y(t)$ denote the number of times state $y$ of arm $i$ is observed up to and including the $t$-th play of arm $i$.
\begin{eqnarray} 
&& P(\bar{r}^i(w_i) \geq \mu^i+c_{t,w_i}) \nonumber \\
&=& P\left( \sum _{y \in S^i} r^i_y n^i_y(w_i) \geq w_i  \sum _{y \in S^i} r^i_y \pi^i_y+w_i c_{t,w_i} \right) \nonumber \\
&=& P\left( \sum_{y \in S^i} (r^i_y n^i_y(w_i) -w_i r^i_y  \pi^i_y ) \geq w_i c_{t,w_i} \right) \nonumber 
\end{eqnarray}
\begin{eqnarray}
\vspace{0.02in}
&=& P\left( \sum_{y \in S^i} (-r^i_y n^i_y(w_i) + w_i r^i_y  \pi^i_y ) \leq - w_i c_{t,w_i} \right) ~. \label{correction1}
\end{eqnarray}
Consider a sample path $\omega$ and the events
\begin{eqnarray*}
A &=&\left\{\omega : \sum_{y \in S^i} (-r^i_y n^i_y(w_i)(\omega) + w_i r^i_y  \pi^i_y ) \leq - w_i c_{t,w_i} \right\} ~,\\
B&=& \bigcup_{y \in S^i} \left\{\omega : -r^i_y n^i_y(w_i)(\omega) + w_i r^i_y  \pi^i_y \leq - \frac{w_i c_{t,w_i}}{|S^i|}\right\} ~. 
\end{eqnarray*}
If $\omega \notin B$, then
\begin{eqnarray*}
&& -r^i_y n^i_y(w_i)(\omega) + w_i r^i_y  \pi^i_y > - \frac{w_i c_{t,w_i}}{|S^i|}, \ \forall y \in S^i \\
\Rightarrow && \sum_{y \in S^i} (-r^i_y n^i_y(w_i)(\omega) + w_i r^i_y  \pi^i_y ) > - w_i c_{t,w_i} ~. 
\end{eqnarray*}
Thus $\omega \notin A$, therefore $P(A) \leq P(B)$. Then continuing from (\ref{correction1}):  
\begin{eqnarray}
&& P(\bar{r}^i(w_i) \geq \mu^i+c_{t,w_i}) \nonumber \\
&\leq& \sum_{y \in S^i} P\left( -r^i_y n^i_y(w_i) + w_i r^i_y  \pi^i_y  \leq - \frac{w_ic_{t,w_i}}{|S^i|} \right) \nonumber
\end{eqnarray}
\begin{eqnarray}
&=& \sum_{y \in S^i} P\left( r^i_y n^i_y(w_i) - w_i r^i_y  \pi^i_y  \geq \frac{w_ic_{t,w_i}}{|S^i|} \right) \nonumber \\
&=& P\left( n^i_y(w_i) -w_i \pi^i_y  \geq \frac{w_ic_{t,w_i}}{|S^i| r^i_y} \right) \nonumber
\end{eqnarray}
\begin{eqnarray}
&=& P\left(\frac{\sum_{t=1}^{w_i} I(X^i_t=y) -w_i \pi^i_y}{\hat{\pi}^i_y w_i} \geq \frac{c_{t,w_i}}{|S^i| r^i_y \hat{\pi}^i_y}\right) \nonumber \\ 
&\leq& \sum_{y \in S^i} N_{q^i} t^{-\frac{L \epsilon^i}{28 (|S^i| r^i_y \hat{\pi}^i_y)^2}}  \label{n1eqn} \\
&\leq& \frac{|S^i|}{\pi_{\min}} t^{-\frac{L \epsilon_{\min}}{28 S^2_{\max} r^2_{\max} \hat{\pi}^2_{\max}}} \label{eqn8},
\end{eqnarray}
where (\ref{n1eqn}) follows from Lemma \ref{lemma3} by letting
\begin{eqnarray*}
\gamma = \frac{c_{t,w_i}}{|S^i| r^i_y \hat{\pi}^i_y}, ~~~ 
f(X^i_t) = \frac{I(X^i_t =y)-\pi^i_y}{\hat{\pi}^i_y} ~, 
\end{eqnarray*}
\rev{and recalling $\hat{\pi}^i_y = \max\{\pi^i_y, 1-\pi^i_y\}$ } 
(note $\hat{P}^i$ is irreducible).


Similarly, we have 
\begin{eqnarray} 
&& P\left(\bar{r}^j(w)\leq\mu^j-c_{t,w}\right) \nonumber \\
&= & P\left(\displaystyle \sum_{y \in S^j} r^j_y (n^j_y(w)-w \pi^j_y) \leq -wc_{t,w}\right) \nonumber \\
&\leq& \sum_{y \in S^j} P\left(r^j_y n^j_y(w)- w r^j_y  \pi^j_y \leq -\frac{w c_{t,w}}{|S^j|} \right) \nonumber \\
&=& \hspace{-4pt}\sum_{y \in S^j}\hspace{-4pt} P\left(r^j_y(w-\displaystyle\sum_{x \neq y} n^j_x(w)) - w r^j_y(1- \displaystyle\sum_{x \neq y} \pi^j_x) \leq -\frac{w c_{t,w}}{|S^j|}\right) \nonumber \\
&=& \sum_{y \in S^j} P\left(r^j_y \displaystyle\sum_{x \neq y} n^j_x(w)-w r^j_y \displaystyle\sum_{x \neq y} \pi^j_x \geq \frac{w c_{t,w}}{|S^j|}\right) \nonumber \\
&\leq& \sum_{y \in S^j} N_{q^j} t^{-\frac{L \epsilon^j}{28 (|S^j|r^j_y \hat{\pi}^j_y)^2}} \label{eqn6} \\
&\leq& \frac{|S^j|}{\pi_{\min}} t^{-\frac{L \epsilon_{\min}}{28 S^2_{\max} r^2_{\max} \hat{\pi}^2_{\max}}} \label{eqn7}  
\end{eqnarray}
where (\ref{eqn6}) again follows from Lemma \ref{lemma3}.
The result \rev{then} follows from combining (\ref{eqn8}) and (\ref{eqn7}): 
\begin{eqnarray}
E\left[\sum_{t=1}^{n} \sum_{w=1}^{t-1} \sum_{w_i=l}^{t-1} I(g^{j}_{t, w} \leq g^i_{t, w_i})\right] &\leq& \frac{|S^i|+|S^j|}{\pi_{\min}} \sum_{t=1}^{\infty}\sum_{w=1}^{t-1}\sum_{w_i=1}^{t-1} t^{-\frac{L \epsilon_{\min}}{28 S^2_{\max} r^2_{\max} \hat{\pi}^2_{\max}}} \nonumber \\
&=&\frac{|S^i|+|S^j|}{\pi_{\min}} \sum_{t=1}^{\infty} t^{-\frac{L \epsilon_{\min} - 56 S^2_{\max} r^2_{\max} \hat{\pi}^2_{\max}}{28 S^2_{\max} r^2_{\max} \hat{\pi}^2_{\max}}} \nonumber \\
&\leq& \frac{|S^i|+|S^j|}{\pi_{\min}} \sum_{t=1}^{\infty} t^{-2}. \label{newlabel1}
\end{eqnarray}
\end{proof}

\section{Proof of Lemma \ref{lemma:rev1}} \label{app:C}

Let $l$ be any positive integer and consider a suboptimal arm $i$. Then,
\begin{eqnarray}
&& T^i(n) = M+ \sum_{t=K+1}^{n} I(i \in A(t)) \leq  M-1+l+\displaystyle\sum_{t=K+1}^{n} I(i \in A(t), T^i(t-1) \geq l)~. \label{eqn:mainrested2}
\end{eqnarray}
Consider
\begin{eqnarray}
E=\bigcup_{j=1}^M \left\{ g^{j}_{t, T^{j}(t)} \leq g^i_{t, T^i(t)} \right\} , \nonumber
\end{eqnarray}
and
\begin{eqnarray}
E^C=\bigcap_{j=1}^M \left\{ g^{j}_{t, T^{j}(t)} > g^i_{t, T^i(t)} \right\}. \nonumber
\end{eqnarray}
If $w \in E^C$ then $i \notin A(t)$. Therefore $\left\{i \in A(t)\right\} \subset E$ and
\begin{eqnarray}
I(i \in A(t), ~T^i(t-1) \geq l) &\leq& I(\omega \in E, ~T^i(t-1) \geq l) \nonumber \\ 
&\leq& \sum_{j=1}^M I(g^{j}_{t, T^{j}(t)} \leq g^i_{t, T^i(t)}, ~T^i(t-1) \geq l) . \nonumber
\end{eqnarray}
Therefore continuing from (\ref{eqn:mainrested2}),
\begin{eqnarray}
T^i(n) &\leq&  M-1+l+ \sum_{j=1}^M \sum_{t=K+1}^{n} I(g^{j}_{t, T^{j}(t)} \leq g^i_{t, T^i(t)}, T^i(t-1) \geq l) \nonumber \\
&\leq&  M-1+l+\sum_{j=1}^M \sum_{t=K+1}^{n} I\left(\min_{1 \leq w \leq t} g^{j}_{t, w} \leq \max_{l \leq w_i \leq t} g^i_{t, w_i} \right) \nonumber \\
&\leq&  M-1+l+ \sum_{j=1}^M \sum_{t=K+1}^{n} \sum_{w=1}^{t-1} \sum_{w_i=l}^{t-1} I(g^{j}_{t,w}\leq g^i_{t, w_i})  \nonumber \\
&\leq& M-1+l + \sum_{j=1}^M \sum_{t=1}^{n} \sum_{w=1}^{t-1} \sum_{w_i=l}^{t-1} I(g^{j}_{t, w} \leq g^i_{t, w_i}) . \nonumber \\
\end{eqnarray} 
Using Lemma \ref{lemma-key} with $l=\left\lceil \frac{4L\ln n}{(\mu^{M}-\mu^i)^2}\right\rceil$, we have for any suboptimal arm
\begin{eqnarray}
E[T^i(n)] \leq M + \frac{4L\ln n}{(\mu^{M}-\mu^i)^2} + \sum_{j=1}^M \frac{(|S^i|+|S^{j}|)\beta}{\pi_{\min}} .\label{eqn:mainrested3}
\end{eqnarray}

\section{} \label{app:D}

\begin{lemma} \label{lemma:analogue}
\rev{Assume Condition \ref{cond:1} holds and all arms are restless.} 
Let $g^i_{t,w}=\bar{r}^i(w)+c_{t,w}$, $c_{t,w}=\sqrt{L\ln t/w}$.  \rev{Under} RCA-M with constant $L \geq 112 S^2_{\max} r^2_{\max} \hat{\pi}^2_{\max} /\epsilon_{\min}$, for any suboptimal arm $i$ and optimal arm $j$ we have 
\begin{eqnarray}
E \left[\sum_{t=1}^{t_2(b)} \sum_{w=1}^{t-1} \sum_{w_i=l}^{t-1} I(g^{j}_{t, w} \leq g^i_{t, w_i}) \right] \leq \frac{|S^i|+|S^j|}{\pi_{\min}} \beta, \label{eqn:new2}
\end{eqnarray}
where \rev{$l=\left\lceil \frac{4L\ln n}{(\mu^M-\mu^i)^2}\right\rceil$ and},
$\beta = \sum_{t=1}^{\infty} t^{-2}$.
\end{lemma}

\begin{proof}
Note that all the quantities in computing the indices in (\ref{eqn:new2}) comes from the intervals $X^i_2(1), X^i_2(2), \cdots \forall i \in \left\{1,\cdots,K\right\}$. Since these intervals begin with state $\gamma^i$ and end with a return to $\gamma^i$ \rev{(but excluding the return visit to $\gamma^i$)}, by the strong Markov property the process at these stopping times have the same distribution as the original process. Moreover by connecting these intervals together we form a continuous sample path which can be viewed as a sample path generated by a Markov chain with an transition matrix identical to the original arm. Therefore we can \rev{proceed in} exactly the same way as the proof of Lemma \ref{lemma-key}. If we choose $s_i \geq 4L \ln(n)/(\mu^M-\mu^i)^2$, then for $t \leq t_2(b)=n' \leq n$, \rev{and} for any suboptimal arm $i$ and optimal arm $j$, 
\begin{eqnarray*}
2c_{t,s_i}=2 \sqrt{\frac{L \ln(t)}{s_i}} \leq 2 \sqrt{\frac{L \ln(t) (\mu^M- \mu^i)^2}{4L\ln(n)}} \leq \mu^j -\mu^i.
\end{eqnarray*}
\rev{The} result follows from letting $l=\left\lceil \frac{4L\ln n}{(\mu^{M}-\mu^i)^2}\right\rceil$ and using Lemma \ref{lemma-key}.
\end{proof}

\section{Proof of Lemma \ref{lemma:restless1}} \label{app:E}

%
%
Let $c_{t,w}=\sqrt{L\ln t/w}$, and let $l$ be any positive integer. Then,
\begin{eqnarray}
&& B^i(b) = 1+\displaystyle\sum_{m=K+1}^{b} I(\alpha(m)=i) \leq  l+\displaystyle\sum_{m=K+1}^{b} I(\alpha(m)=i, B^i(m-1) \geq l)  \label{transactions1}
\end{eqnarray}
Consider any sample path $\omega$ and the following sets
\begin{eqnarray}
E=\bigcup_{j=1}^M \left\{\omega : g^{j}_{t_2(m-1), T^{j}_2(t_2(m-1))}(\omega) \leq g^i_{t_2(m-1), T^i_2(t_2(m-1))}(\omega) \right\}, \nonumber
\end{eqnarray}
and
\begin{eqnarray}
E^C=\bigcap_{j=1}^M \left\{\omega : g^{j}_{t_2(m-1), T^{j}_2(t_2(m-1))}(\omega) > g^i_{t_2(m-1), T^i_2(t_2(m-1))}(\omega) \right\}. \nonumber
\end{eqnarray}
If $\omega \in E^C$ then $\alpha(m) \neq i$. Therefore $\left\{\omega: \alpha(m)(\omega)=i\right\} \subset E$ and
\begin{eqnarray}
&& I(\alpha(m)=i, B^i(m-1) \geq l) \leq I(\omega \in E, B^i(m-1) \geq l) \nonumber \\ 
&\leq& \sum_{j=1}^M I(g^{j}_{t_2(m-1), T^{j}_2(t_2(m-1))} \leq g^i_{t_2(m-1), T^i_2(t_2(m-1))}, B^i(m-1) \geq l)~. \nonumber
\end{eqnarray}
Therefore continuing from (\ref{transactions1}),
\begin{eqnarray}
B^i(b) &\leq&  l+ \sum_{j=1}^M \sum_{m=K+1}^{b} I(g^{j}_{t_2(m-1), T^{j}_2(t_2(m-1))} \leq g^i_{t_2(m-1), T^i_2(t_2(m-1))}, B^i(m-1) \geq l) \nonumber \\
&\leq&  l+\sum_{j=1}^M \sum_{m=K+1}^{b} I\left(\min_{1 \leq w \leq t_2(m-1)} g^{j}_{t_2(m-1), w} \leq \max_{t_2(l) \leq w_i \leq t_2(m-1)} g^i_{t_2(m-1), w_i} \right) \nonumber \\
&\leq&  l+ \sum_{j=1}^M \sum_{m=K+1}^{b} \sum_{w=1}^{t_2(m-1)} \sum_{w_i=t_2(l)}^{t_2(m-1)} I(g^{j}_{t_2(m),w}\leq g^i_{t_2(m), w_i})  \label{jul21} \\
&\leq& l + M \sum_{j=1}^M \sum_{t=1}^{t_2(b)} \sum_{w=1}^{t-1} \sum_{w_i=l}^{t-1} I(g^{j}_{t, w} \leq g^i_{t, w_i})~, \label{eqn1} 
\end{eqnarray} 
where as given in (\ref{eqn1redifine}),  
$g^i_{t,w}=\bar{r}^i(w)+c_{t,w}$, \rev{and we have assumed that the index value of an arm remains the same between two updates.} 
The inequality in (\ref{eqn1}) follows from the facts that the second outer sum in (\ref{eqn1}) is over time while the second outer sum in (\ref{jul21}) is over blocks, each block lasts at least two time slots and at most $M$ blocks can be completed \rev{in} each time step. From this point on we use Lemma \ref{lemma:analogue} to get
\begin{eqnarray}
E[B^i(b(n))|b(n)=b] \leq \left\lceil \frac{4L\ln t_2(b)}{(\mu^{M}-\mu^i)^2}\right\rceil + M \sum_{j=1}^M \frac{(|S^i|+|S^{j}|)\beta}{\pi_{\min}},\nonumber
\end{eqnarray}
for all suboptimal arms. Therefore,
\begin{eqnarray}
E[B^i(b(n))] \leq \frac{4L\ln n}{(\mu^{M}-\mu^i)^2} +1+ M \sum_{j=1}^M C_{i,j} \beta, \label{eqn9}
\end{eqnarray}
since $n \geq t_2(b(n))$ almost surely.


The total number of plays of arm $i$ at the end of block $b(n)$ is equal to the total number of plays of arm $i$ during the regenerative cycles of visiting state $\gamma^i$ plus the total number of plays before entering the regenerative cycles plus one more play resulting from the last play of the block which is state $\gamma^i$.  This gives: 
\begin{eqnarray}
E[T^i(n)] \leq \left(\frac{1}{\pi^i_{\min}} + \Omega^i_{\max} +1\right)E[B^i(b(n))] ~. \nonumber
\end{eqnarray}
Thus,
\begin{eqnarray}
&& \sum_{i>M} (\mu^1-\mu^i) E[T^i(n)] \\
& \leq&  4L \sum_{i>M} \frac{(\mu^1-\mu^i) D_i \ln n}{(\mu^M-\mu^i)^2} + \sum_{i>M} (\mu^1-\mu^i) D_i \left(1 + M \sum_{j=1}^M   C_{i,j}\right). \label{eqn10}
\end{eqnarray}

\section{Proof of Theorem \ref{thm:restless1}} \label{app:F}

Assume that the states which determine the regenerative sample paths are given {\em a priori} by $\gamma=[\gamma^1,\cdots,\gamma^K]$. \rev{This is to simplify the analysis by skipping the initialization stage of the algorithm and we will show that this choice does not affect the regret bound.} We denote the expectations with respect to \rev{RCA-M} given $\gamma$ as $E_\gamma$. First we rewrite the regret in the following form:
\begin{eqnarray}
R_\gamma(n) &=& \sum_{j=1}^M \mu^{j} E_\gamma[T(n)] - E_\gamma \left[ \sum_{t=1}^{T(n)} \sum_{\alpha(t) \in A(t)} r^{{\alpha}(t)}_{{x}_{{\alpha}(t)}}\right] 
+ \sum_{j=1}^M \mu^{j} E_\gamma[n-T(n)]-E_\gamma \left[\sum_{t=T(n)+1}^{n} \sum_{\alpha(t) \in A(t)} r^{\alpha(t)}_{x_{\alpha(t)}}\right] \nonumber \\
&=& \left\{ \sum_{j=1}^M \mu^{j} E_\gamma[T(n)] - \sum_{i=1}^K \mu^i E_\gamma \left[T^i(n)\right]  \right\} - Z_{\gamma}(n) \label{eqn:last_R} \\
&+& \sum_{j=1}^M \mu^{j} E_\gamma[n-T(n)]-E_\gamma \left[\sum_{t=T(n)+1}^{n} \sum_{\alpha(t) \in A(t)} r^{\alpha(t)}_{x_{\alpha(t)}}\right], \label{eqn:last_R1}
\end{eqnarray}
where for notational convenience, we have used
\begin{eqnarray}
Z_\gamma(n)=E_\gamma \left[\sum_{t=1}^{T(n)} \sum_{\alpha(t) \in A(t)}  r^{\alpha(t)}_{x_{\alpha(t)}} \right]  - \sum_{i=1}^K \mu^i E_\gamma \left[T^i(n)\right]. \nonumber
\end{eqnarray}
We have
\begin{eqnarray}
 \sum_{j=1}^M \mu^{j} E_\gamma[T(n)] - \sum_{i=1}^K \mu^i E_\gamma \left[T^i(n)\right] 
 &=& \sum_{j=1}^M \sum_{i=1}^K \mu^{j} E_\gamma[T^{i,j}(n)] - \sum_{j=1}^{\rev{M}} \sum_{i=1}^K \mu^i E_\gamma[T^{i,j}(n)] \nonumber \\
&=& \sum_{j=1}^M \sum_{i>M} (\mu^{j}-\mu^i) E_\gamma [T^{i,j}(n)] \nonumber \\
&\leq& \sum_{i>M} (\mu^{1}-\mu^i) E_\gamma [T^i(n)] \label{transactions2}
\end{eqnarray}

Since we can bound (\ref{transactions2}), i.e. the difference in \rev{the} brackets in (\ref{eqn:last_R}) logarithmically using Lemma \ref{lemma:restless1}, it remains to bound $Z_\gamma(n)$ and the difference in (\ref{eqn:last_R1}). We have
\begin{eqnarray}
Z_\gamma(n) &\geq& \sum_{i=1}^M \sum_{y \in S^{i}}r^{i}_y E_\gamma\left[\sum_{b=1}^{B^{i}(b(n))} \sum_{X^{i}_t \in X^{i}(b)} I(X^{i}_t=y)\right] \nonumber \\ 
&& + \sum_{i>M} \sum_{y \in S^i} r^i_y E_\gamma \left[ \sum_{b=1}^{B^i(b(n))} \sum_{X^i_t \in X^i_2(b)} I(X^i_t=y)\right] \label{modif1} \\
&& - \sum_{i=1}^M \mu^{i} E_\gamma \left[T^{i}(n)\right] \nonumber \\
&& - \sum_{i>M} \mu^i \left(\frac{1}{\pi^i_{\gamma^i}}+\Omega^i_{\max}+1\right) E_\gamma \left[ B^i(b(n)) \right] ~, \nonumber
\end{eqnarray}
where the inequality comes from counting only the rewards obtained during the \rev{SB2's} for all suboptimal arms and the last part of the proof of Lemma \ref{lemma:restless1}. 
Applying Lemma \ref{lastlemma} to (\ref{modif1}) we get 
\begin{eqnarray}
E_\gamma \left[ \sum_{b=1}^{B^i(b(n))} \sum_{X^i_t \in X^i_2(b)} I(X^i_t=y)\right] = \frac{\pi^i_y}{\pi^i_{\gamma^i}} E_\gamma \left[B^i(b(n))\right] ~.  \nonumber
\end{eqnarray}
Rearranging terms we get 
\begin{eqnarray}
Z_\gamma(n) \geq R^*(n) - \sum_{i>M} \mu^i (\Omega^i_{\max}+1) E_\gamma \left[ B^i(b(n)) \right] \label{modif3}
\end{eqnarray}
where
\begin{eqnarray}
R^*(n)&=& \sum_{i=1}^M \sum_{y \in S^{i}} r^{i}_y E_\gamma \left[ \sum_{b=1}^{B^{i}(b(n))} \sum_{X^{i}_t \in X^{i}(b)}  I(X^{i}_t=y)\right] 
- \sum_{i=1}^M \sum_{y \in S^{i}} r^{i}_y \pi^{i}_y  E_\gamma \left[T^{i}(n)\right] \nonumber.
\end{eqnarray}

Consider now $R^*(n)$. Since all suboptimal arms are played at most logarithmically, the total number of time slots in which an optimal arm is not played is at most logarithmic.  It follows that the number of discontinuities between plays of any single optimal arm is at most logarithmic. 
For any optimal arm $i \in \left\{1,\cdots,M\right\}$ we combine \rev{{\em consecutive}} blocks in which \rev{arm $i$} is played into a single \rev{{\em combined}} block, and denote by $\bar{X}^{i}(j)$ the $j$-th combined block of \rev{arm $i$}. 
Let $\bar{b}^{i}$ denote the total number of combined blocks for arm $i$ up to block $b$. Each $\bar{X}^{i}$ thus consists of two sub-blocks: $\bar{X}^{i}_1$ that contains the states visited from \rev{the} beginning of $\bar{X}^{i}$ (empty if the first state is $\gamma^{i}$) to the state right before hitting $\gamma^{i}$, and sub-block $\bar{X}^{i}_2$ that contains the rest of $\bar{X}^{i}$ (a random number of regenerative cycles).  

\comment{
\begin{figure}
\includegraphics[width=3.5in]{figureb.eps}
\caption{$\bar{X}^*, \bar{X}^*_1, \bar{X}^*_2$ for a sample path}
\label{figure5}
\end{figure}
} 

Since a \rev{combined} block $\bar{X}^{i}$ \rev{necessarily} starts after \rev{certain} discontinuity in playing the \rev{$i$-th} best arm, $\bar{b}^{i}(n)$ is less than or equal to the total number of \rev{discontinuities} of play \rev{of} the \rev{$i$-th} best arm up to time $n$.  \rev{At the same time,} the total number of \rev{discontinuities} of play \rev{of} the \rev{$i$-th} best arm up to time $n$ is less than or equal to the total number of blocks in which suboptimal arms are played up to time $n$.  Thus 
\begin{eqnarray}
E_\gamma[\bar{b}^{i}(n)] \leq \sum_{k>M} E_\gamma[B^k(b(n))]. \label{modif4}
\end{eqnarray}
We \rev{now} rewrite  $R^*(n)$ in the following from:
\begin{eqnarray}
R^*(n)&=& \sum_{i=1}^M \sum_{y \in S^{i}} r^{i}_y E_\gamma \left[ \sum_{b=1}^{\bar{b}^{i}(n)} \sum_{X^{i}_t \in \bar{X}^{i}_2(b)}  I(X^{i}_t=y)\right] \label{difference6} \\
&&- \sum_{i=1}^M \sum_{y \in S^{i}} r^{i}_y \pi^{i}_y  E_\gamma \left[ \sum_{b=1}^{\bar{b}^{i}(n)} |\bar{X}^{i}_2(b)| \right] \label{difference7} \\
&&+ \sum_{i=1}^M \sum_{y \in S^{i}} r^{i}_y E_\gamma \left[ \sum_{b=1}^{\bar{b}^{i}(n)} \sum_{X^{i}_t \in \bar{X}^{i}_1(b)}  I(X^{i}_t=y)\right] \label{modif5} \\
&&- \sum_{i=1}^M \sum_{y \in S^{i}} r^{i}_y \pi^{i}_y  E_\gamma \left[ \sum_{b=1}^{\bar{b}^{i}(n)} |\bar{X}^{i}_1(b)| \right] \label{modif6} \\
&&> 0 - \sum_{i=1}^M \mu^{i} \Omega^{i}_{\max} \sum_{k>M} E_\gamma[B^k(b(n))] \label{difference8}
\end{eqnarray}
where the last inequality is obtained by noting the difference between (\ref{difference6}) and (\ref{difference7}) is zero by Lemma \ref{lastlemma}, using positivity of rewards to lower bound (\ref{modif5}) by $0$, and (\ref{modif4}) to upper bound (\ref{modif6}). 
Combining this with (\ref{eqn9}) and (\ref{modif3}) we can obtain a logarithmic upper bound on $-Z_\gamma(n)$ by the following steps:
\begin{eqnarray}
-Z_\gamma(n) &\leq& -R^*(n) + \sum_{i>M} \mu^i (\Omega^i_{\max}+1) E_\gamma \left[ B^i(b(n)) \right] \nonumber \\
&\leq& \sum_{i=1}^M \mu^{i} \Omega^{i}_{\max} \sum_{k>M} \left(\frac{4L\ln n}{(\mu^{M}-\mu^k)^2} +1+ M \sum_{j=1}^M C_{k,j} \beta \right) \nonumber \\
&+& \sum_{i>M} \mu^i (\Omega^i_{\max}+1) \left(\frac{4L\ln n}{(\mu^{M}-\mu^i)^2} +1+ M \sum_{j=1}^M C_{k,i} \beta \right) \nonumber \\
\end{eqnarray}  
We also have,
\begin{eqnarray}
\sum_{j=1}^M \mu^{j} E_\gamma[n-T(n)]-E_\gamma[\sum_{t=T(n)+1}^{n} \sum_{\alpha(t) \in A(t)} r^{\alpha(t)}_{x_{\alpha(t)}}]
&\leq& \sum_{j=1}^M \mu^{j} E_\gamma[n-T(n)] \nonumber \\
&=& \sum_{j=1}^M \mu^{j}\left(\frac{1}{\pi_{\min}} + \max_{i \in \left\{1,...,K\right\}} \Omega^i_{\max} +1\right). \label{last_eq}
\end{eqnarray}

\rev{Finally, combining the above results as well as Lemma \ref{lemma:restless1} we get}
\begin{eqnarray}
R_\gamma(n) &=& \left\{ \sum_{j=1}^M \mu^{j} E_\gamma[T(n)] - \sum_{i=1}^K \mu^i E_\gamma \left[T^i(n)\right]  \right\} - Z_{\gamma}(n) \nonumber \\
&&+ \sum_{j=1}^M \mu^{j} E_\gamma[n-T(n)]-E_\gamma \left[\sum_{t=T(n)+1}^{n} \sum_{\alpha(t) \in A(t)} r^{\alpha(t)}_{x_{\alpha(t)}}\right] \nonumber \\
&\leq& \sum_{i>M} (\mu^{1}-\mu^i) E_\gamma [T^i(n)] \nonumber \\
&&+ \sum_{i=1}^M \mu^{i} \Omega^{i}_{\max} \sum_{k>M} \left(\frac{4L\ln n}{(\mu^{M}-\mu^k)^2} +1+ M \sum_{j=1}^M C_{k,j} \beta \right) \nonumber \\
&&+ \sum_{i>M} \mu^i (\Omega^i_{\max}+1) \left(\frac{4L\ln n}{(\mu^{M}-\mu^i)^2} +1+ M \sum_{j=1}^M C_{k,i} \beta \right) \nonumber \\
&&+  \sum_{j=1}^M \mu^{j}\left(\frac{1}{\pi_{\min}} + \max_{i \in \left\{1,...,K\right\}} \Omega^i_{\max} +1\right) \nonumber \\
&=& 4L \ln n \sum_{i>M} \frac{1}{(\mu^{M}-\mu^i)^2}\left((\mu^{1}-\mu^i) D_i + E_i\right) \nonumber \\
&&+ \sum_{i>M} \left((\mu^{1}-\mu^i)D_i+E_i\right) \left(1+ M \sum_{j=1}^M C_{i,j}\right) +F \nonumber
\end{eqnarray}
Therefore we have obtained the stated logarithmic bound for (\ref{eqn:last_R}).  Note that this bound does not depend on $\gamma$, and therefore is also an upper bound for $R(n)$, completing the proof. 


\end{document}